\documentstyle[txmac,a4,
amssymb,%
case,%
twoside,%
nocaphead,%
epsf,%
varthm,%
rotate,%
myrot,%
mypic,times,mathptm]{article}

\advance\oddsidemargin by -0.7cm
\advance\evensidemargin by -0.7cm
\advance\textwidth by 1.4cm

\def\mynewtheo#1#2{%
\newtheorem{@#1}{#2}[section]%
\newenvironment{#1}{\begin{@#1}\rm}{\end{@#1}}}

\mynewtheo{lemma}{Lemma}
\mynewtheo{exer}{Exercise}
\mynewtheo{theo}{Theorem}
\mynewtheo{rem}{Remark}
\mynewtheo{defi}{Definition}
\mynewtheo{conj}{Conjecture}
\mynewtheo{corr}{Corollary}
\mynewtheo{prop}{Proposition}
\mynewtheo{ques}{Question}
\mynewtheo{exam}{Example}

\newenvironment{eqn}{\begin{equation}}{\end{equation}}

\begin{document}

\makeatletter

\newenvironment{myeqn*}[1]{\begingroup\def\@eqnnum{\reset@font\rm#1}%
\xdef\@tempk{\arabic{equation}}\begin{equation}\edef\@currentlabel{#1}}
{\end{equation}\endgroup\setcounter{equation}{\@tempk}\ignorespaces}

\newenvironment{myeqn}[1]{\begingroup\let\eq@num\@eqnnum
\def\@eqnnum{\bgroup\let\r@fn\normalcolor 
\def\normalcolor####1(####2){\r@fn####1#1}%
\eq@num\egroup}%
\xdef\@tempk{\arabic{equation}}\begin{equation}\edef\@currentlabel{#1}}
{\end{equation}\endgroup\setcounter{equation}{\@tempk}\ignorespaces}

\newcommand{\mybin}[2]{\text{$\Bigl(\begin{array}{@{}c@{}}#1\\#2%
\end{array}\Bigr)$}}
\newcommand{\mybinn}[2]{\text{$\biggl(\begin{array}{@{}c@{}}%
#1\\#2\end{array}\biggr)$}}

\def\overtwo#1{\mbox{\small$\mybin{#1}{2}$}}
\newcommand{\mybr}[2]{\text{$\Bigl\lfloor\mbox{%
\small$\displaystyle\frac{#1}{#2}$}\Bigr\rfloor$}}
\def\mybrtwo#1{\mbox{\mybr{#1}{2}}}

\def\myfrac#1#2{\raisebox{0.2em}{\small$#1$}\!/\!\raisebox{-0.2em}{\small$#2$}}

\def\myeqnlabel{\bgroup\@ifnextchar[{\@maketheeq}{\immediate
\stepcounter{equation}\@myeqnlabel}}

\def\@maketheeq[#1]{\def\theequation{#1}\@myeqnlabel}

\def\@myeqnlabel#1{%
{\edef\@currentlabel{\theequation}
\label{#1}\enspace\eqref{#1}}\egroup}

\def\epsfs#1#2{{\epsfxsize#1\relax\epsffile{#2.eps}}}


\author{A. Stoimenow\\[2mm]
\small Humboldt University Berlin, Dept.~of Mathematics,\\
\small Ziegelstra\ss e 13a, 10099 Berlin, Germany,\\
\small e-mail: {\tt stoimeno@informatik.hu-berlin.de},\\
\small WWW: {\hbox{\tt http://www.informatik.hu-berlin.de/\raisebox{-0.8ex}{\tt\~{}}stoimeno}}
}

\title{\large\bf \uppercase{Positive knots, closed braids}\\[2mm]
\uppercase{and the Jones polynomial}}

\date{\large Current version: \today\ \ \ First version:
\makedate{5}{5}{1997}}



\maketitle

\makeatletter

\def\chrd#1#2{\picline{1 #1 polar}{1 #2 polar}}
\def\gchrd#1#2{\picline{1.5 #1 polar}{1.5 #2 polar}}
\def\arrow#1#2{\picvecline{1 #1 polar}{1 #2 polar}}

\def\labch#1#2#3{\chrd{#1}{#2}\picputtext{1.3 #2 polar}{$#3$}}
\def\labar#1#2#3{\arrow{#1}{#2}\picputtext{1.3 #2 polar}{$#3$}}
\def\labbr#1#2#3{\arrow{#1}{#2}\picputtext{1.3 #1 polar}{$#3$}}

\def\labline#1#2#3#4{\picvecline{#1}{#2}\pictranslate{#2}{
  \picputtext{#3}{$#4$}}}
\def\lablineb#1#2#3#4{\picvecline{#1}{#2}\pictranslate{#1}{
  \picputtext{#3}{$#4$}}}

\def\pt#1{{\picfillgraycol{0}\picfilledcircle{#1}{0.06}{}}}
\def\labpt#1#2#3{\pictranslate{#1}{\pt{0 0}\picputtext{#2}{$#3$}}}

\def\CD#1{{\let\@nomath\@gobble\small\diag{6mm}{2}{2}{
  \picveclength{0.4}\picvecwidth{0.14}
  \pictranslate{1 1}{
    \piccircle{0 0}{1}{}
    #1
}}}}

\def\GD#1{{\let\@nomath\@gobble\scriptsize\diag{3mm}{3.0}{3.0}{
  \picveclength{0.32}\picvecwidth{0.14}
  \pictranslate{1.5 1.5}{
    \piccircle{0 0}{1}{}
    #1
}}}}

\def\conf#1{
\let\@conf\thr@@conf
\ifcase#1\def\@array{[0 60 300 {}f]}
\or\def\@array{[0 70 110 180 250 290]}
\or\def\@array{[50 90 130 230 270 310]}
\or\def\@array{[0 90 180 270]}
\let\@conf\tw@conf\fi\@conf}
\def\thr@@conf#1#2#3#4#5#6#7#8#9{\GD{\ea\picPSgraphics\ea{\@array /@A x D}%
\@mkchr{#1}{#2}{#7}
\@mkchr{#3}{#4}{#8}
\@mkchr{#5}{#6}{#9}
}}
\def\tw@conf#1#2#3#4#5#6{\GD{\ea\picPSgraphics\ea{\@array /@A x D}%
\@mkchr{#1}{#2}{#5}
\@mkchr{#3}{#4}{#6}
}}
\def\@mkchr#1#2#3{\ifnum#1<0\let\@tempy\labch\else\let\@tempy\labar
\fi\@tempy{@A #1 abs 1 - g}{@A #2 1 - g}{#3}}

\let\fa\forall
\let\ay\asymp
\let\pa\partial
\let\al\alpha
\let\be\beta
\let\bt\beta
\let\Gm\Gamma
\let\gm\gamma
\let\de\delta
\let\dl\delta
\let\Dl\Delta
\let\eps\epsilon
\let\lm\lambda
\let\Lm\Lambda
\let\sg\sigma
\let\vp\varphi
\let\om\omega

\let\sm\setminus
\let\tl\tilde
\def\ncap{\not\mathrel{\cap}}
\def\disc{\text{\rm disc}\,}
\def\cf{\text{\rm cf}\,}
\def\sp{\text{\rm span}\,}
\def\lra{\longrightarrow}
\def\so{\Rightarrow}
\def\So{\Longrightarrow}
\let\ds\displaystyle

\let\reference\ref

\long\def\@makecaption#1#2{%
   \vskip 10pt
   {\let\label\@gobble
   \let\ignorespaces\@empty
   \xdef\@tempt{#2}%
   }%
   \ea\@ifempty\ea{\@tempt}{%
   \setbox\@tempboxa\hbox{%
      \fignr#1#2}%
      }{%
   \setbox\@tempboxa\hbox{%
      {\fignr#1:}\capt\ #2}%
      }%
   \ifdim \wd\@tempboxa >\captionwidth {%
      \rightskip=\@captionmargin\leftskip=\@captionmargin
      \unhbox\@tempboxa\par}%
   \else
      \hbox to\captionwidth{\hfil\box\@tempboxa\hfil}%
   \fi}%
\def\fignr{\small\sffamily\bfseries}%
\def\capt{\small\sffamily}%

\newdimen\@captionmargin\@captionmargin2cm\relax
\newdimen\captionwidth\captionwidth\hsize\relax

\def\eqref#1{(\protect\ref{#1})}

\def\proof{\@ifnextchar[{\@proof}{\@proof[\unskip]}}
\def\@proof[#1]{\noindent{\bf Proof #1.}\enspace}

\def\hint{\noindent Hint: }
\def\problem{\noindent{\bf Problem.} }
\def\note{\noindent{\bf Note.} }
\def\question{\noindent{\bf Question.} }

\def\@mt#1{\ifmmode#1\else$#1$\fi}
\def\qed{\hfill\@mt{\Box}}
\def\qqed{\hfill\@mt{\Box\enspace\Box}}

\def\cU{{\cal U}}
\def\cC{{\cal C}}
\def\cP{{\cal P}}
\def\tP{{\tilde P}}
\def\tZ{{\tilde Z}}
\def\fg{{\frak g}}
\def\tr{\text{tr}}
\def\cZ{{\cal Z}}
\def\cV{{\cal V}}
\def\cD{{\cal D}}
\def\bR{{\Bbb R}}
\def\cE{{\cal E}}
\def\bZ{{\Bbb Z}}
\def\bN{{\Bbb N}}

\def\bysame{\same[\kern2cm]\,}

\def\br#1{\left\lfloor#1\right\rfloor}
\def\BR#1{\left\lceil#1\right\rceil}

\def\abstractname{}

\@addtoreset {footnote}{page}

\renewcommand{\section}{%
   \@startsection
         {section}{1}{\z@}{-1.5ex \@plus -1ex \@minus -.2ex}%
               {1ex \@plus.2ex}{\large\bf}%
}
\renewcommand{\@seccntformat}[1]{\csname the#1\endcsname .
\quad}

\def\bC{{\Bbb C}}
\def\bP{{\Bbb P}}

{\let\@noitemerr\relax
\vskip-2.7em\kern0pt\begin{abstract}
\noindent{\bf Abstract.}\enspace
Using the recent Gauss diagram formulas for Vassiliev invariants
of Polyak-Viro-Fiedler and combining these formulas
with the Bennequin inequality, we prove several inequalities for
positive knots relating their Vassiliev invariants, genus
and degrees of the Jones polynomial. As a consequence,
we prove that for any of the polynomials of Alexander/Conway,
Jones, HOMFLY, Brandt-Lickorish-Millett-Ho and
Kauffman there are only finitely many positive knots
with the same polynomial and no positive knot with trivial polynomial.

We also discuss an extension
of the Bennequin inequality, showing that the  unknotting number
of a positive  knot not less than its genus, which recovers
some recent unknotting number results of A'Campo, Kawamura and Tanaka,
and give applications to the Jones polynomial of a positive knot.
\\[1mm]
\noindent\em{Keywords:} positive knots, Vassiliev invariants,
Gau\ss{} sums, Jones polynomial, Alexander polynomial, genus,
Casson invariant, unknotting number.\\[1mm]
\end{abstract}
}

{\parskip0.2mm\tableofcontents}
\vspace{7mm}

\section{\label{sect1}Introduction}

Positive knots, the knots having diagrams with all crossings positive,
have been for a while of interest for knot theorists, not only because
of their intuitive defining property. Such knots have occurred, in the
more special case of braid positive knots (in this paper knots
which are closed positive braids will be called so) in the theory
of dynamical systems \cite{BirWil}, singularity theory \cite{ACampo,%
BoiWeb}, and in the more general class of quasipositive knots (see
\cite{Rudolph3}) in the theory of algebraic curves \cite{Rudolph}.

Beside the study of some classical invariants of positive knots
\cite{Busk,CochranGompf,Traczyk}, significant progress in the study
of such knots was achieved by the discovery of the new polynomial
invariants \cite{Jones,HOMFLY,Kauffman2,BLM,Ho}, giving rise to
a series of results on properties of these invariants for this knot
class \cite{Cromwell,Fiedler,MorCro,Yokota,Zulli}.

Recently, a conceptually new approach for defining invariants
of finite type (Vassiliev invariants) \cite{BirmanLin,BarNatanVI,%
BarNatanBibl,BarNatanOde,survey,Bseq,Vogel,Vassiliev} was initiated
by Fiedler \cite{Fiedler,Fiedler2,Fiedler3} and Polyak-Viro
\cite{VirPol} by the theory of small state (or Gau\ss{}) sums.
Fiedler remarked \cite{Fiedler4} that the Gau\ss{} sum formulas
have direct application to the study of positive knots. 

This paper aims to work out a detailed account on such applications.
Sharpening Fiedler's results, we will prove a number of inequalities
for positive knots, relating via the Gau\ss{} sum formulas the
Vassiliev invariants $v_2(K)$ and $v_3(K)$ of degree $2$ and $3$ of
a positive knot $K$ on the one hand, and classical invariants like
its genus $g(K)$, crossing number $c(K)$ and unknotting number $u(K)$
on the other hand. We will use the tables of Rolfsen \cite[appendix]
{Rolfsen} and Thistlethwaite \cite{KnotScape} to find examples 
illustrating and showing
the essence of these properties as positivity criteria for knots.

Although all inequalities can be considered in their own right, one
of them, which subsequently turned out of central importance, and
is thus worth singling out, is the inequality $v_2(K)\ge c(K)/4$ we
will prove in \S\reference{secCas}. Similar (although
harder to prove) inequalities will be first discussed in
\S\reference{iv3} for $v_3$, improving the one originally
given by Fiedler. As a consequence of involving the crossing
number into our bounds, we prove, that there are only finitely many
positive knots with the same Jones polynomial $V$ and that any knot has
only finitely many (possibly no) positive reduced diagrams, so that
positivity can always be (at least theoretically) decided, provided
one can identify a knot from a given diagram. A further application of
such type of inequality is given in \cite{gen1}, where it is decisively
used to give polynomial bounds of the number of positive knots
of fixed genus and given crossing number.

For our results on unknotting numbers it will also turn out useful
to apply the machinery of inequalities of Bennequin type \cite[theorem
3, p.~101]{Bennequin} for the (slice) genus. Thus we devote
a separate section \S\reference{un} to the discussion of this topic.
In particular, there we give an extension of Bennequin's inequality
to arbitrary diagrams. An application of this
extension is the observation that the unknotting number
of a positive knot is not less than its genus (corollary
\reference{corr2}). This resolves, \em{inter alia}, the unknotting
numbers of $5$ of the undecided knots in Kawauchi's tables
\cite{Kawauchi}, which have been (partially) obtained by
Tanaka \cite{Tanaka}, Kawamura \cite{Kawamura} and A'Campo
\cite{ACampo} (examples \reference{_ex1}, \reference{_ex2} and
\reference{_ex3}). It can also be used to extend some results proved
on the genus of positive knots to their unknotting number (see
\cite{gen1}).


%
%


In \S\reference{rp} we will use the inequalities derived for
$v_2$ and $v_3$ together with those given by Morton \cite{Morton2}
to give some relations between the values of $v_2$, $v_3$ and
the HOMFLY polynomial of positive knots.

Braid positive knots, \em{inter alia}, because of their special
importance will be considered in their own right in \S\reference{S.3},
where some further specific inequalities for the Vassiliev invariants
will be given.  We will also prove, that the minimal degree of the
Jones polynomial of a closed positive braid is equal for knots
to the genus and is at least a quarter of its crossing number.

Finally, in the sections \reference{un} and \reference{qu} 
we will review some results and conjectures and summarize some
questions, which are interesting within our setting.

{\bf Notation.}
For a knot $K$ denote by $c(K)$ its (minimal) crossing number, by $g(K)$
its genus, by $b(K)$ its braid index, by $u(K)$ its unknotting number,
by $\sg(K)$ its signature.
$!K$ denotes the obverse (mirror image) of $K$. We use the
Alexander-Briggs notation and the Rolfsen \cite{Rolfsen} tables
to distinguish between a knot and its obverse. ``Projection''
is the same as ``diagram'', and this means a knot or link
diagram. Diagrams are always assumed oriented.

The symbol $\Box$ denotes the end or the absence of a proof.
In latter case it is \em{assumed} to be evident from the
preceding discussion/references; else (and anyway)
I'm grateful for any feedback.

\section{Positive knots and Gau\ss{} sums}

\begin{defi}
The writhe is a number ($\pm1$), assigned to any crossing in a link
diagram. A crossing as on figure \ref{figwr}(a), has writhe 1 and
is called positive. A crossing as on figure \ref{figwr}(b), has writhe 
$-1$ and is called negative. A crossing is smoothed out by replacing
it by the fragment on figure \ref{figwr}(c) (which changes the number
of components of the link). A crossing as on figure \reference{figwr}%
(a) and \reference{figwr}(b) is smashed to a singularity (double
point) by replacing it by the fragment on figure \reference{figwr}(d).
A $m$-singular diagram is a diagram with $m$ crossings smashed.
A $m$-singular knot is an immersion prepresented by a
$m$-singular diagram.

\end{defi}

\begin{figure}[htb]
\[
\begin{array}{c@{\qquad}c@{\qquad}c@{\qquad}c}
\diag{6mm}{1}{1}{
\picmultivecline{0.12 1 -1.0 0}{1 0}{0 1}
\picmultivecline{0.12 1 -1.0 0}{0 0}{1 1}
} &
\diag{6mm}{1}{1}{
\picmultivecline{0.12 1 -1 0}{0 0}{1 1}
\picmultivecline{0.12 1 -1 0}{1 0}{0 1}
} &
\diag{6mm}{1}{1}{
\piccirclevecarc{1.35 0.5}{0.7}{-230 -130}
\piccirclevecarc{-0.35 0.5}{0.7}{310 50}
} &
\diag{6mm}{1}{1}{
\picmultivecline{0.12 1 -1.0 0}{1 0}{0 1}
\picmultivecline{0.12 1 -1.0 0}{0 0}{1 1}
\picfillgraycol{0}
\picfilledcircle{0.5 0.5}{0.1}{}
}
\\[2mm]
(a) & (b) & (c) & (d)
\end{array}
\]
\caption{\label{figwr}}
\end{figure}

\begin{defi}
A knot is called positive, if it has a positive diagram,
i.~e. a diagram with all crossings positive.
\end{defi}

Recall \cite{NewInv,VirPol} the concept of Gau\ss{} sum invariants.
As they will be the main tool of all the further investigations,
we summarize for the benefit of the reader the basic points of
this theory.

\begin{figure}[htb]
{
\[
\begin{array}{c@{\qquad}c}
     \diag{2cm}{2.0}{2.0}{
  \picputtext[dl]{0 0}{\epsfs{4cm}{6_2}}
  \picputtext[u]{1.25 1.65}{$1$} 
  \picputtext[u]{0.8 1.65}{$2$} 
  \picputtext[u]{0.6 0.6}{$3$} 
  \picputtext[u]{1.0 0.35}{$4$} 
  \picputtext[u]{1.3 0.9}{$5$} 
  \picputtext[u]{0.9 1.3}{$6$} 
	}
 & \input{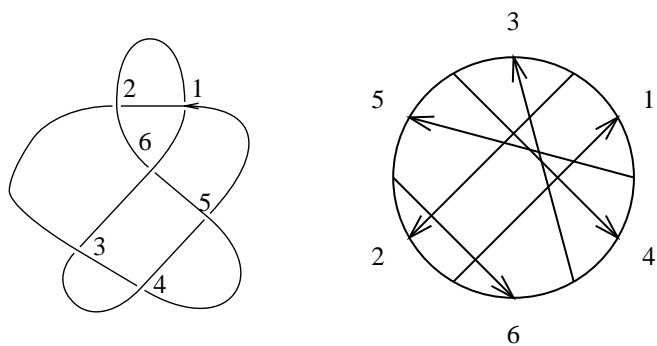}
\end{array}
\]}
\caption{The knot $6_2$ and its Gau\ss{} diagram.\label{fig6_2}}
\end{figure}

\begin{@defi}[\protect\cite{Fiedler3,VirPol}]\rm
A Gau\ss{} diagram (GD) of a knot diagram is an oriented circle with
arrows connecting points on it mapped to a crossing and
oriented from the preimage of the undercrossing to the
preimage of the overcrossing. See figure \reference{fig6_2}.
\end{@defi}

Fiedler \cite{Fiedler3,NewInv} found the following formula for (a
variation of) the degree-3-Vassiliev invariant using Gau\ss{} sums.
\begin{eqn}\label{v3}
v_3\,=\,\sum_{(3,3)}w_pw_qw_r\,+\,
\sum_{(4,2)0}w_pw_qw_r\,+\,\frac{1}{2}\sum_\text{$p,q$ linked}(w_p+w_q)\,,
\end{eqn}
where the configurations are
\[
\begin{array}{c@{\qquad}c@{\qquad}c}
\CD{\chrd{50}{240}
\chrd{90}{270}
\chrd{125}{320}} &
\CD{\chrd{-180}{0}
\arrow{-70}{70}
\arrow{-110}{110}} &
\CD{\labar{250}{70}{$q$}
\labar{290}{110}{$p$}} \\[4mm]
(3,3) & (4,2)0 & \text{a linked pair}
\end{array}
\]
Here chords depict arrows which may point in both directions and
$w_p$ denotes the writhe of the crossing $p$. For a given
configuration, the summation in \eqref{v3} is done over each
unordered pair/triple of crossings, whose arrows in the
Gau\ss{} diagram form that configuration. The terms associated to a
pair/triple of crossings occurring in the sums are called \em{weights}.
If no weight is specified, we take by default the product of the
writhes of the involved crossings. Thus
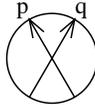
\[
\CD{\labar{240}{60}{$q$}
\labar{300}{120}{$p$}}
\]
means `sum of $w_p\cdot w_q$ over $p,q$ linked'. In the linked pair
of the picture above, call $p$ \em{distinguished}, that is,
the over-crossing of $p$ is followed by the under-crossing of $q$.
For the motivation of this notation, see \cite{NewInv}.

Additionally, one may put a base point on both the knot
and Gau\ss{} diagram (see \cite{VirPol}). This is equivalent to
distinguishing a cyclic order of the arrow ends, or ``cutting''
the circle somewhere.

To make precise which variation of the degree-3-Vassiliev invariant
we mean, we noted in \cite{NewInv}, that
\[
v_3\,=\,-\frac 13V^{(2)}(1)-\frac 19V^{(3)}(1)\,,
\]
where $V$ is the Jones polynomial \cite{Jones} and $V^{(n)}$
denotes the $n$-th derivation of $V$.
We noted further (and shall use it later), that $v_3$ is additive under
connected knot sum, that is, $v_3(K_1\#K_2)=v_3(K_1)+v_3(K_2)$ (verify
this!).

\begin{defi} A diagram is composite, if it looks as in figure
\ref{figtan}(a) and both $A$ and $B$ contain at least one crossing.
A diagram is split, if it looks as in figure
\ref{figtan}(b) and both $A$ and $B$ are non-empty.
A composite link is a link with a composite diagram, in which no one
of $A$ and $B$ represent the unknot. A split link is a link with a 
split diagram.
\end{defi}

We will use the synonyms `prime' and `connected' for
`non-composite' and `disconnected' for `composite'.

\begin{figure}[htb]
\[
\begin{array}{*2c}
\diag{6mm}{4}{2}{
  \piccirclearc{2 0.5}{1.3}{45 135}
  \piccirclearc{2 1.5}{1.3}{-135 -45}
  \picfilledcircle{1 1}{0.8}{A}
  \picfilledcircle{3 1}{0.8}{B}
} &
\diag{6mm}{4}{2}{
  \picfilledcircle{1 1}{0.8}{A}
  \picfilledcircle{3 1}{0.8}{B}
} \\
(a) & (b)
\end{array}
\]
\caption{\label{figtan}}
\end{figure}

\begin{defi} A crossing is reducible, if its smoothing out
yields a split diagram. A diagram is reduced, if it has no
reducible crossings.
\end{defi}

\begin{defi}
Call a positive diagram bireduced, if it is reduced and does not
admit a move
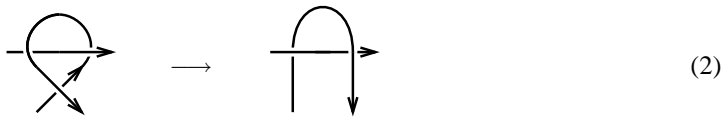
\begin{eqn}\label{second}
\diag{1cm}{2}{1.6}{
  \picline{0.3 1}{1 1}
  \picmulticurve{0.12 1 -1.0 0}{1.0 1.5}{0.7 1.5}{0.4 1.1}{0.7 0.8}
  \picmulticurve{0.12 1 -1.0 0}{1.3 0.8}{1.6 1.1}{1.3 1.5}{1.0 1.5}
  \picmultivecline{0.12 1 -1.0 0}{1 1}{1.7 1}
  \picvecline{0.7 0.2}{1.3 0.8}
  \picmultivecline{0.12 1 -1.0 0}{0.7 0.8}{1.3 0.2}
  \piccurve{1.0 1.5}{0.7 1.5}{0.4 1.1}{0.7 0.8}
}\quad\lra\quad
\diag{1cm}{2}{1.6}{
  \picvecline{1 1}{1.7 1}
  \picmulticurve{0.12 1 -1.0 0}{1.4 1}{1.4 1.8}{0.6 1.8}{0.6 1}
  \picmultivecline{0.12 1 -1.0 0}{1.4 1}{1.4 0.2}
  \picline{0.6 0.2}{0.6 1}
  \picmultiline{0.12 1 -1.0 0}{0.3 1}{1 1}
  \picline{0.9 1}{1.1 1}
}
\end{eqn}
To this move we will henceforth refer as a second (reduction) move.
\end{defi}

The reason for introducing this move will become clear shortly.

\begin{defi}
The intersection graph of a Gau\ss{} diagram is a graph with vertices corresponding to
arrows in the Gau\ss{} diagram and edges connecting intersecting arrows/vertices.
\end{defi}

Gau\ss{} diagrams have in general the following properties.

\begin{defi}
For two chords in a Gau\ss{} diagram $a\cap b$ means
``$a$ intersects $b$'' (or crossings $a$ and $b$ are linked)
and $a\ncap b$ means ``$a$ does not intersect $b$''
(or crossings $a$ and $b$ are not linked).
\end{defi}

We now formulate two simple properties of Gau\ss{} diagrams
that will be extensively used in the following, even valence
and double connectivity.

\begin{@lemma}[double connectivity]
\label{lem1}\rm
Whenever in a Gau\ss{} diagram $a\cap c$ and $b\cap c$
then either $a\cap b$ or there is an arrow $d$ with
$d\cap a$ and $d\cap b$. I.~e., in the intersection graph
of the Gau\ss{} diagram
any two neighbored edges participate in a cycle of length 3 or 4.
In particular, the Gau\ss{} diagram (or its intersection graph)
are doubly connected.
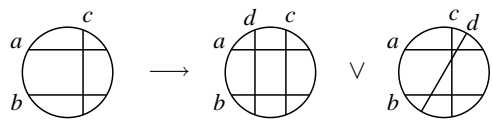
\[
\CD{\labch{-70}{70}{c}
\labch{30}{150}{a}
\labch{-30}{-150}{b}
} \quad\lra\quad
\CD{\labch{-70}{70}{c}
\labch{30}{150}{a}
\labch{-30}{-150}{b}
\labch{-110}{110}{d}
} \quad\lor\quad
\CD{\labch{-80}{80}{c}
\labch{30}{150}{a}
\labch{-30}{-150}{b}
\labch{240}{60}{d}
}
\]
\end{@lemma}

\proof Assume $a\ncap b$. Consider the plane curve of the projection.
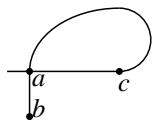
\[
\diag{6mm}{4}{3}{
  \picveclength{0.4}\picvecwidth{0.14}
  \picline{0.5 1}{3 1}
  \piccirclearc{3 1.7}{0.7}{-90 90}
  \picellipsearc{3 1}{2 1.4}{90 180}
  \picline{1 1}{1 0}
  \labpt{1 1}{0.2 -0.2}{a}
  \labpt{3 1}{0.1 -0.3}{c}
  \labpt{1 0}{0.2 0.2}{b}
}
\]

As the curve meets $c$ the second time before doing so with $b$,
it has a segment in the inner part of
the above depicted loop between both occurrences of $a$,
and so there must be another crossing between the first and second
occurence of  $a$ and the first and second occurence of  $b$. \qed

\section{Inequalities for $v_3$\label{iv3}}

In this section we shall prove an obstruction to
positivity which renders it decidable, whether a given knot
has this property. The idea is due to Fiedler, but here we present
an improved version of it.

Our goal is now to prove the following two statements.
%
{\nopagebreak
\def\labelenumi{\theenumi)}\mbox{}\\[-18pt]

\begin{enumerate}
\item\label{item1} The number of edges in the intersection graph of a non-composite
Gau\ss{} diagram (=intersections of chords in the Gau\ss{} diagram=linked pairs) is at least
$3\mybrtwo{c-1}$, where $c$ is the number of vertices in the intersection graph
(=chords in the Gau\ss{} diagram=crossings in the knot projection).
\item\label{item2} In any positive diagram $D$ of $c$ crossings,
$v_3(D)\ge \#\{\text{ linked pairs }\}\ge c$.
If $D$ is bireduced, then $v_3(D)\ge \myfrac{4}{3}\#\{\text{
linked pairs }\}$.
%
\end{enumerate}
}

We do this in steps and split the arguments into several lemmas.
Finally, we summarize the results in a more self-contained form
in theorem \ref{th1}. We start by

\begin{lemma}\label{lp1}
If $K$ is a positive reduced diagram of $c$ crossings,
then $v_3(K)\ge c$.
\end{lemma}

\proof Consider the Gau\ss{} diagram of this projection.

By the additivity of $v_3$ and $c$ in composite projections
(note, that factors in composite positive reduced projections
and themselves in positive reduced projections),
henceforth assume, the positive diagram is non-composite,
i.~e. the Gau\ss{} diagram of $c$ arrows is connected.

But then this diagram will have at least $c-1$ intersections, i.~e.
linked pairs, and the above argument (lemma \reference{lem1}) shows that, as the number
of arrows is more than $2$, the number of chord intersections
cannot be equal to $c-1$. As in a positive diagram each linked pair
contributes one to the value of the third term in \eqref{v3} and the first
two terms are non-negative, the assertion follows. \qed

Bounds of this kind render it decidable whether a positive diagram
exists.

\begin{corr}
Any reduced positive diagram of any knot $K$ has maximally $v_3(K)$
crossings. In particular, there are only finitely many positive
knots with the same $v_3$ and any knot has only finitely many 
(possibly no) positive diagrams. \qed
\end{corr}

Therefore, there are also only finitely many positive knots with 
the same Jones polynomial, but we will state this fact in greater
generality somewhat later.

\begin{corr}
If $K$ is not the unknot, maximally one of $K$ and $!K$
can be positive. In particular, no positive non-trivial knot
is amphicheral.
\end{corr}

This fact follows in the special case of alternating knots from
Thistlethwaite's invariance of the writhe \cite{Kauffman2}
and for Lorenz knots from work of Birman and Williams \cite{BirWil}.
A general proof was first given by Cochran and Gompf
\cite[corollary 3.4, p.~497]{CochranGompf} and briefly later
independently by Traczyk \cite{Traczyk} using the signature. 
We will later, in passing by, give an independent argument for
the positivity of the signature on positive knots using Gau\ss{}
diagrams.

\proof It follows from $v_3(!K)=-v_3(K)$ which is easy to see from the
formula \eqref{v3}: mirroring reverses the orientation
of all arrows in the Gau\ss{} diagram and all configurations in \eqref{v3}
are invariant under this operation, while the terms in the sum change
the sign. \qed

\begin{exam} $4_1$ (the figure eight knot) and $6_3$ are
amphicheral and hence cannot be positive.
\end{exam}

\begin{rem}
The bound of lemma \reference{lp1} is sharp, as $v_3(!3_1)=4$ and there is the following
reduced positive 4 crossing diagram  of the right-hand trefoil:
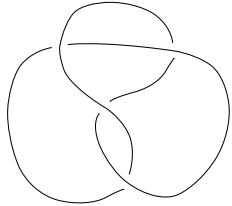
\[
\epsfs{3cm}{k_tref1}
\]
However, this diagram is not bireduced and here I came
to consider this notion.
\end{rem}

\begin{exam}
Beside the standard and the above depicted diagram, there cannot
be more positive reduced diagrams of the right-hand trefoil $!3_1$.
\end{exam}

\begin{exam}
Following T.~Fiedler, and as indicated in \cite{NewInv},
the knot $6_2$ has $v_3=4$. So it cannot have any positive diagram
(as else it would have a reduced one and this would have to have not
more than $4$ crossings).
\end{exam}

Here is another property of Gau\ss{} diagrams we will use
in the following to sharpen our bound. 

\begin{@lemma}[even valence]\label{lem2} \rm
Any chord in a Gau\ss{} diagram has odd length (i.~e., even number of
basepoints on both its sides, or equivalently,
even number of intersections with other chords, that is, even valence
in the intersection graph of the Gau\ss{} diagram).
\end{@lemma}

\proof This is, as lemma \ref{lem1}, a consequence of
the Jordan curve theorem, and is reflected e.~g. also
in the definition of the Dowker notation of knot diagrams
\cite{ThiDow}. \qed


Here is the improved bound announced in
\cite{NewInv} under assumption of bireducedness.

\begin{lemma}
If $K$ is a positive bireduced diagram of $c$ crossings, then
\[
v_3(K)\ge \myfrac{4}{3}\,\#\text{ linked pairs }\ge
\myfrac{4}{3}\,c\,.
\]
In particular, $v_3(K)\ge \myfrac{4}{3}c(K)$ for $K$ positive.
\end{lemma}

\proof To prove is the first inequality (the second was proved in lemma
\reference{lp1}). Assume w.l.o.g. as before the
Gau\ss{} diagram is connected. We know that the
number of intersections in the Gau\ss{} diagram ($=$ number of linked pairs)  is at least
$c$. So it suffices to prove 
\[
\#\{\text{ matching $(3,3)$ and $(4,2)0$ configurations }\}\,
\ge \myfrac{1}{3}\#\{\text{ linked pairs }\}\,.
\]
To do this, we will construct a map
\[
m\,:\,\{\text{ crossings in the GD (linked pairs) }\}\,\lra
\,\{\text{ matching $(3,3)$ and $(4,2)0$ configurations }\}
\]
such that each image is realized not more than $3$ times. To prove this
property of $m$, we will check it each time we define a new
value of $m$ on the values of $m$ defined so far.

About the definition of $m$. Set $m$ on a crossing
participating in a $(3,3)$ or $(4,2)0$ configuration
to one (any arbitrary)  of these configurations.
So, up to now, all $(4,2)0$ configurations are realized as
image under $m$ maximally $2$ and all $(3,3)$ configurations
are realized as image under $m$ maximally $3$ times.

Now look at a crossing $A$, not participating in \em{any}
$(3,3)$ and $(4,2)0$ configuration.
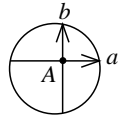
\[
\CD{
\labar{170}{10}{a}
\labar{-80}{80}{b}
\labpt{0.18 0.18}{-0.3 -0.3}{A}
}
\]
If chord $a$ has length $3$ then we have either 

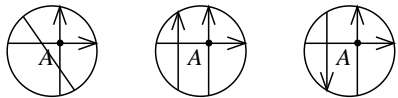
\[
{
\def\CDA#1{\CD{\arrow{170}{10}\arrow{-80}{80}
\labpt{0.18 0.18}{-0.3 -0.3}{A}#1}}
%
\CDA{\chrd{130}{-60}} \qquad
\CDA{\arrow{-120}{120}} \qquad
\CDA{\arrow{120}{-120}}
}
\]
In the first two cases $A$ is in a $(3,3)$ or $(4,2)0$ configuration,
and in the third case this is exactly the situation of a second move
\eqref{second}.
Note: it follows from the positivity of the diagram, that indeed
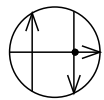
\[
\CD{\arrow{-180}{0}\arrow{65}{-65}\arrow{-120}{120}\pt{0.45 0}
\pt{0.45 0}}
\]
does not exist. Else the diagram part on the left in \eqref{second}
to be positive, we had to reverse the direction of (exactly)
one of the strands, and the crossings would become linked.
So let $a$ have length at least $5$.

\begin{caselist}

\def\CDA#1{
  \CD{
    \labch{-110}{110}{a}
    \labch{-160}{-20}{b}
    \labch{150}{30}{y}
    \labpt{-0.33 -0.33}{-0.25 -0.25}{A}
    #1
}}
\case \label{case0} First assume $a$ has only $2$ crossings.
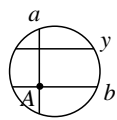
\[
\CDA{}
\]
We have $y\ncap b$ (else $A\in (3,3)$). By double connectivity
$\exists\, x, x\cap y,b$.
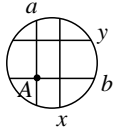
\[
\CDA{\labch{80}{-80}{x}}
\]
$x$ does not intersect $a$ (else $A\in (3,3)$).

On the other hand, if for some $c$, $c\cap y$, then $c\cap b$ and vice versa
(else by double connectivity on $a,c,y$ we had
$\exists d\cap y, d\cap a$. As $b\ncap y$ we had $d\neq b$ and so
$d$ would be a third intersection of $a$).
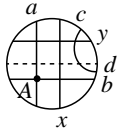
\[
\CDA{\labch{80}{-80}{x}
  \piccurve{1 50 polar}{0.3 0.3}{0.6 -0.2}{1 -10 polar}
  \picputtext{1  50 polar 0.3 +}{$c$}
  {\piclinedash{0.1}{0.05}\labch{-180}{0}{d}}
}
\]
As $a$ is not of length $3$, on the other side of $a$
from that, where $x$ lies, there must be a chord $z$ which
(by assumption of connectedness of the diagram) must intersect
one of $b$ or $y$ and therefore (see above) both.
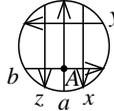
\[
\CD{\labar{65}{-65}{x}
\labar{115}{-115}{z}
\labbr{-90}{90}{a}
\labbr{30}{150}{y}
\labbr{-150}{-30}{b}
\labpt{0 -0.5}{0.2 -0.2}{A}
}
\]
Then $z,x$ must be equally oriented with respect to $y$ and
$b$ (else $A\in (4,2)0$), i.~e. $(z,x,y)$ and $(z,b,x)$
are of type $(4,2)0$. Assign by $m$ to $A$ the second one of these
configurations. \myeqnlabel{*1}

\case So now let $a$ have at least $4$ crossings (remember, each chord has
even number of crossings!).
Look at $a$:
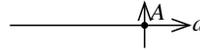
\[
  \diag{6mm}{4}{1}{\small
    \picveclength{0.4}\picvecwidth{0.14}
    \labline{0 0.5}{4 0.5}{0.2 0}{a}
    \picvecline{3 0}{3 1}
    \labpt{3 0.5}{0.3 0.3}{A}
  }
\]
Beside by $b$, $a$ is intersected $n\ge 3$ times by (only) downward pointing
arrows (else either $A\in (3,3)$ or $A\in (4,2)0$).

\begin{caselist}

\case \label{case1}%
Two such chords $a_1,a_2$ do not intersect.
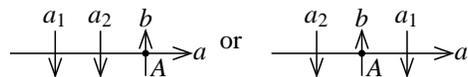
\[
  \diag{6mm}{4}{1.5}{
    \picveclength{0.4}\picvecwidth{0.14}
    \labline{0 0.5}{4 0.5}{0.2 0}{a}
    \labline{3 0}{3 1}{0 0.3}{b}
    \lablineb{2 1}{2 0}{0 0.3}{a_2}
    \lablineb{1 1}{1 0}{0 0.3}{a_1}
    \labpt{3 0.5}{0.3 -0.3}{A}
  }\quad
  \text{or}\quad
  \diag{6mm}{4}{1.5}{
    \picveclength{0.4}\picvecwidth{0.14}
    \labline{0 0.5}{4 0.5}{0.2 0}{a}
    \labline{2 0}{2 1}{0 0.3}{b}
    \lablineb{3 1}{3 0}{0 0.3}{a_1}
    \lablineb{1 1}{1 0}{0 0.3}{a_2}
    \labpt{2 0.5}{0.3 -0.3}{A}
  }
\]
Set
\begin{eqn}\label{*2}
m(A)=K:=\{a_1,a_2,a\}\in (4,2)0\,.
\end{eqn}
Up to now, $K\in (4,2)0$ has only 2 preimages, unless it was not
the object of an assignment of the kind \eqref{*1} or \eqref{*2}
before. However, there is only maximally one such additional
preimage $A$ of $K$, because we can uniquely reconstruct $A$ from
$K$:
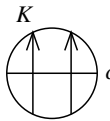
\[
\CD{\arrow{-115}{115}
\arrow{-65}{65}
\labch{-180}{0}{c}
\picputtext{-0.6 1.3}{$K$}
}
\]
Consider the chord $c$ in $K$ with both intersections on it.
Then the other two arrows point in 1 direction with respect to
$c$. $A$ is then the unique intersection point of $c$ with an
arrow pointing in the opposite direction than the other two arrows
of $K$ do.

Summarizing cases \ref{case0} and \ref{case1},
no $(4,2)0$ configuration received more than $3$ preimages
as far as $m$ is constructed now.

\case \label{case2}
All $n\ge 3$ chords intersect. The picture is like this
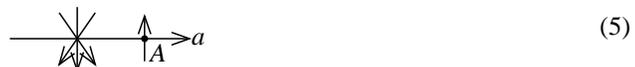
\begin{eqn}\label{five}
\diag{6mm}{4}{1.5}{
  \picveclength{0.4}\picvecwidth{0.14}
  \labline{0 0.5}{4 0.5}{0.2 0}{a}
  \picvecline{3 0}{3 1}
  \labpt{3 0.5}{0.3 -0.3}{A}
  \pictranslate{1.5 0.5}{
    \picmultigraphics[rt]{3}{35}{
      \picvecline{0.7 55 polar}{0.7 235 polar}
    }
  }
}
\end{eqn}
These $n$ chords produce with $a$ $\overtwo{n}$ configurations
of type $(3,3)$ and among themselves $\overtwo{n}$ intersections.
So $n+\overtwo{n}$ intersection points participate in $\overtwo{n}$
configurations $(3,3)$ involving $a$. So there is a relation among
these $\overtwo{n}$ with a preimage under $m$ of maximally
\[
\frac{n+\overtwo{n}}{\overtwo{n}}\le 2
\]
intersections participating in the configuration.
Define by $m$ on $A$ as any of these relations.

\end{caselist}
\end{caselist}

How many preimages now has a configuration of type $(3,3)$?
If it was not affected by the so far considered configurations
in case \ref{case2}, it still has maximally $3$ preimages.
If it has been, it has maximally $2$ preimages among the intersections
participating in it. How many ``$A$''s could have been assigned
to such a configuration $K$ by case \ref{case2}? If any, 
$K$ must look like
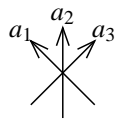
\[
\diag{6mm}{2}{2}{
  \picveclength{0.4}\picvecwidth{0.14}
  \pictranslate{1 1}{
  \labar{225}{45}{a_3}
  \labar{270}{90}{a_2}
  \labar{315}{135}{a_1}
}}
\]
and $A$ must be either on $a_3$ or $a_1$ and be the unique intersection
point of a chord intersecting $a_1$ (resp.~$a_3$) in the reverse
direction as all other chords, among others, $a_3$ (resp.~$a_1$),
do (as this chord is different from $a_3$ (resp.~$a_1$),
its intersection direction is uniquely determined).
So there are at most 2 such ``$A$''s and the configuration has
at most 4 preimages.

%

We would like to show now that in fact $(3,3)$ configurations with
4 preimages can always be avoided by a proper choice of
$(3,3)$ configurations in case \reference{case2}.

Assume, that at one point in case \ref{case2} all
configurations $(3,3)$ of $a$ with two downward pointing arrows
in \eqref{five} already have $3$ preimages
as a next $A$ has to be added (that is, you
are forced to create a fourth preimage to one of the
$(3,3)$ configurations). Then there is only one choice.
There are exactly $3$ chords (which mutually intersect and
intersect $a$), from the resulting $6$ crossings and $3$
configurations $(3,3)$ involving $a$, each configuration contains
exactly $2$ of its points in  its preimage (for $n>3$ we have
\[
\frac{n+\overtwo{n}}{\overtwo{n}}<2\,,
\]
and so there is always a configuration with not more than one
of its points in its preimage) and to each of these $3$ configurations
$(3,3)$ there has already been assigned an ``$A$'' by case \ref{case2}.
(Here ``$A$'' means an intersection point, which participated as $A$ in
some previous application of case \ref{case2}.)
There cannot have been 2 ``$A$''s added, as $A$ would be the third
possible one and we saw that there are no 3 possible ones for the same
$(3,3)$ configuration.
Because on each chord of the configuration only one possible
``$A$'' can lie, this other ``$A$'' (different from our $A$)
must lie on
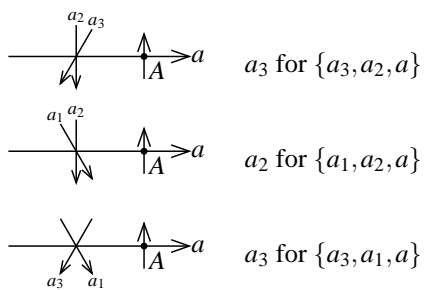
\[
\def\CDA#1{
\diag{6mm}{4}{2.0}{
  \picveclength{0.4}\picvecwidth{0.14}
  \pictranslate{0 0.7}{
  \labline{0 0.5}{4 0.5}{0.2 0}{a}
  \picvecline{3 0}{3 1}
  \labpt{3 0.5}{0.3 -0.3}{A}
  \pictranslate{1.5 0.5}{
    \picscale{0.7 0.7}{
  #1
  }}}
}}
\begin{array}{c@{\qquad}c}
\CDA{ \labbr{60}{240}{a_3}
      \labbr{90}{270}{a_2}
    }
 & a_3 \text{ for } \{a_3,a_2,a\}\\
\CDA{ \labbr{120}{300}{a_1}
      \labbr{90}{270}{a_2}
    }
 & a_2 \text{ for } \{a_1,a_2,a\}\\
\CDA{ \labar{60}{240}{a_3}
      \labar{120}{300}{a_1}
    }
 & a_3 \text{ for } \{a_3,a_1,a\}
\end{array}
\]
But this cannot be, because to the ``$A$'' on $a_3$
(it is unique, because there's always $a$ in the configuration and this
``$A$'' must intersect with $a_3$ in the opposite direction,
and 2 such ``$A$''s would $\in (4,2)0$) cannot simultaneously have
been assigned both $\{a_3,a_1,a\}$ and $\{a_3,a_2,a\}$ under $m$.
This contradiction shows, that it must be really always possible
to define $m$ on an ``$A$'' in case \ref{case2}, not augmenting
the number of preimages of a $(3,3)$ configuration to more than $3$.

So now any configuration of type $(3,3)$ has maximally $3$ preimages and $m$ is
completely defined, and has the desired property. \qed


But of course, there are in general much more linked pairs
than crossings, and so we can go a little further.

Consider the intersection graph $G$ of a Gau\ss{} diagram.

\begin{lemma}
In $G$
\[
\#\text{ edges }\,\ge \,3\,\left(\br{\frac{\#\text{ vertices } -1}{2} }
\right)\,,
\]
if $G$ connected, i.~e. the Gau\ss{} diagram non-composite.
\end{lemma}


\proof Recall that the intersection graph of a Gau\ss{} diagram has the 
double connectivity property, that each pair of neighbored edges
lies in some 3 or 4 cycle.
\begin{eqn}\label{my*}
\diag{6mm}{1}{1}{
  \picline{0 0}{1 0}
  \picline{0 0}{0.3 0.8}
  \labpt{0 0}{-0.2 -0.2}{b}
  \labpt{1 0}{0.2 -0.2}{c}
  \labpt{0.3 0.8}{-0.2 0.2}{a}
} \quad\lra\quad
\diag{6mm}{1}{1}{
  \picline{0 0}{1 0}
  \picline{0 0}{0.3 0.8}
  \labpt{0 0}{-0.2 -0.2}{b}
  \labpt{1 0}{0.2 -0.2}{c}
  \labpt{0.3 0.8}{-0.2 0.2}{a}
  \picline{1 0}{0.3 0.8}
} \quad\text{ or }\quad
\diag{6mm}{1}{1}{
  \picline{0 0}{1 0}
  \picline{0 0}{0 1}
  \picline{1 1}{1 0}
  \picline{1 1}{0 1}
  \labpt{0 0}{-0.2 -0.2}{b}
  \labpt{1 0}{0.2 -0.2}{c}
  \labpt{1 1}{0.2 0.2}{d}
  \labpt{0 1}{-0.2 0.2}{a}
}
\end{eqn}
Fix a spanning tree $B$ of the intersection graph $G$. Denote by $c$
the number of vertices in $G$. $B$ has $c-1$ edges, as $G$ is connected.
Choose a disjoint cover $\cU$ of $\,2\mybrtwo{c-1}$ edges in $B$
in $\mybrtwo{c-1}$ pairs, so that each pair has neighbored edges
(why does this work?). Apply \eqref{my*} to construct for each pair
a map
\[
m\,:\,\cU'\,\lra\,\{\text{ edges and pairs of edges in $G$
outside $B$ }\}\,,
\]
where $\cU'$ is an extension of $\cU$ (i.~e. $\forall A\in\cU\,\exists\,
A'\in\cU':A'\supset A$).
Define $m$ as follows. Fix a pair
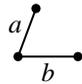
\def\CDA#1{
\diag{8mm}{1}{1}{
  \picline{0 0}{1 0}
  \picline{0 0}{0.3 0.8}
  \pt{0 0}\pt{1 0}\pt{0.3 0.8}
  \picputtext[r]{0.06 0.5}{$a$}
  \picputtext[u]{0.5 -0.1}{$b$}
  #1
}}
\[
\CDA{}
\]
in $B$. There are 3 cases (the full edges belong to $B$ and the
dashed edges are outside of $B$).
\[
1) \quad
\CDA{\picputtext[dl]{0.7 0.4}{$c$}
{\piclinedash{0.1}{0.05}
  \picline{1 0}{0.3 0.8}
}}
  \qquad 2) \quad
\diag{8mm}{1}{1}{
  \picline{0 0}{1 0}
  \picline{0 0}{0 1}
  {\piclinedash{0.1}{0.05}
  \picline{1 1}{1 0}
  \picline{1 1}{0 1}
  }
  \pt{0 0}\pt{1 0}\pt{1 1}\pt{0 1}
  \picputtext[r]{-0.1 0.5}{$a$}
  \picputtext[u]{0.5 -0.1}{$b$}
  \picputtext[l]{1.1 0.5}{$d$}
  \picputtext[d]{0.5 1.1}{$c$}
}  
  \qquad 3) \quad
\diag{8mm}{1}{1}{
  \picline{0 0}{1 0}
  \picline{0 0}{0 1}
  \picline{1 1}{1 0}
  {\piclinedash{0.1}{0.05}
  \picline{1 1}{0 1}
  }
  \pt{0 0}\pt{1 0}\pt{1 1}\pt{0 1}
  \picputtext[r]{-0.1 0.5}{$a$}
  \picputtext[u]{0.5 -0.1}{$b$}
  \picputtext[l]{1.1 0.5}{$d$}
  \picputtext[d]{0.5 1.1}{$c$}
}  
\]
\begin{caselist}

\case\label{_case1} Set $m(\{a,b\}):=c$.

\case\label{_case2} Set $m(\{a,b\}):=\{c,d\}$.

\case\label{_case3} Set $m(\{a,b,d\}):=c$.

\end{caselist}

As all the pairs are disjoint, all triples obtained by extending a pair
by one element in case \reference{_case3} are distinct, and the map will be well defined.
$\cU'$ is the cover $\cU$,
where some of the pairs have been extended to triples by case
\ref{_case3}. We have $\#\,\cU'=\mybrtwo{c-1}$.

\begin{figure}[htb]
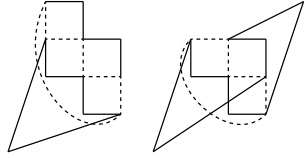

\[
\begin{array}{*2c}
\diag{5mm}{3}{4}{
  \picline{1 3}{1 2}
  \piclineto{2 2}
  \piclineto{2 1}
  \piclineto{3 1}
  \piclineto{0 0}
  \piclineto{1 3}
  \picline{1 4}{2 4}
  \piclineto{2 3}
  \piclineto{3 3}
  \piclineto{3 2}
  \piclinedash{0.1}{0.05}
  \picline{1 4}{1 3}
  \piclineto{2 3}
  \piclineto{2 2}
  \piclineto{3 2}
  \piclineto{3 1}
  \piccurveto{2 0}{0 2}{1 4}
} &
\diag{5mm}{4}{4}{
  \picline{1 3}{1 2}
  \piclineto{2 2}
  \piclineto{2 1}
  \piclineto{3 1}
  \piclineto{4 4}
  \piclineto{2 3}
  \piclineto{3 3}
  \piclineto{3 2}
  \piclineto{0 0}
  \piclineto{1 3}
  \piclinedash{0.1}{0.05}
  \picline{1 3}{2 3}
  \piclineto{2 2}
  \piclineto{3 2}
  \piclineto{3 1}
  \piccurveto{2 0}{0 2}{1 3}
}
\end{array}
\]
\caption{How a cycle in $F$ implies one in $B$. Full edges belong to
$B$ and dashed outside of $B$. The vertices in the cycle in $F$
are the edges in the dashed cycle in $G$. The edges in the cycle in $F$
are pairs of neighbored edges in the dashed cycle in $G$. The
edges in $B$ depicted belong to the pairs of $\cU$, which are preimages
of the edges in the cycle in $F$. They are all disjoined and
contain together a cycle.\label{figFB}}
\end{figure}

As $B$ is without cycle, no edge in $G\sm B$ has received two preimages
by cases \ref{_case1} and \ref{_case3}. In the same way, no pair
of (neighbored) edges in $G\sm B$ received two preimages
by case \ref{_case2}. Moreover, if we look at the dual graph of $G$
(where such pairs correspond to edges), the edges in $G\sm B$ with a
preimage by case \reference{_case2} from a forest $F$.
(Convince yourself, using figure \reference{figFB},
that the existence of a cycle in $F$ implies one
in $B$.) Therefore, for all components $C$ of $F$ the number of
involved vertices in $C$ (=edges in $G\sm B$ involved in one of
these pairs) is bigger that the number of edges of $C$ (=pairs of
edges in $G\sm B$ with a preimage by $m$). Furthermore, for all
components $C$ of $F$ maximally one of the vertices of $C$ (=edges
in $G\sm B$) has a preimage by \ref{_case1}) or \ref{_case3})
(again as $B$ is a tree). So we see that
\[
\#\,\cU'\le \#\{ \text{ edges in $G\sm B$ }\}.
\]
So there are at least $\mybrtwo{c-1}$ edges in $G\sm B$ and at least
\[
c-1 + \mybrtwo{c-1}\ge 3\mybrtwo{c-1}
\]
edges in $G$.  \qed

Summarizing, we proved:


\begin{theo}\label{th1}
If $K$ is a positive bireduced non-composite projection of $c$
crossings, then $v_3(K)\ge 4\mybrtwo{c-1}$. If the projection is
not bireduced, but reduced and non-composite, we have at least $v_3(K)
\ge 3\mybrtwo{c-1}$. If the projection is composite and bireduced,
we have $v_3(K)\ge \frac{4}{3}\,c$, and if it is
composite and reduced but not bireduced, $v_3(K)\ge c$. \qed
\end{theo}

\begin{corr}
If a prime knot $K$ is positive, it has a positive diagram of
not more than $v_3(K)/2+2$ crossings. \qed
\end{corr}

There exist some other generally sharper obstructions to positivity.
One is due to Morton and Cromwell \cite{MorCro}: If $P$
denotes the HOMFLY polynomial \cite{HOMFLY} 
(in the convention of \cite{polynomial}), then for a positive link
$P(it,iz)$ must have only non-negative coefficients
in $z$ for any $t\in[0,1]$ ($i$ denotes $\sqrt{-1}$).
The special case $t=1$ is the positivity
of the Conway polynomial \cite{Conway},
proved previously for braid positive
links by v.~Buskirk \cite{Busk} and later extended to positive
links by Cromwell \cite{Cromwell}.

Moreover, in \cite{Cromwell} it was proved, that for $L$ positive
$\min\deg_l(P)=\max\deg_m(P)$.

That these obstructions, although generally sharper, are not
always better, shows the following example, coming out of some quest 
in Thistlethwaite's tables.

\begin{figure}[htb]
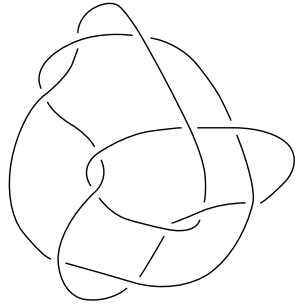

\[
\epsfs{4cm}{k-12-2038}
\]
\caption{\label{fig12-2038}The knot $12_{2038}$.}
\end{figure}

\begin{exam}
The knot $12_{2038}$ on figure \reference{fig12-2038} has the
HOMFLY polynomial
\[
(-7l^6-9l^8-3l^{10})+(13l^6+13l^8+3l^{10})m^2+(-7l^6-6l^8-l^{10})m^4
+(l^6+l^8)m^6\,.
\]
It shows, that the obstructions of \cite{Busk,Cromwell} and
\cite{MorCro} are not violated. However, $v_3(12_{2038})=8$.
\end{exam}

\begin{rem}
One may ask, in how far can the given bounds be improved.
The answer is, using our arguments, not very much, as shows
the following
\end{rem}

\begin{exam}
Consider the graph $G_n$, which is the Hasse diagram of the lattice
$(\cP(\{1,$ $\dots,n\}),\subset)$. I.~e., its vertices are subsets
of $\{1,\dots,n\}$ and $A$ and $B$ are connected by an edge, if
$B\subset A$ and $\#(A\sm B)=1$.
Then $G_n$ satisfies the double connectivity property
of lemma \ref{lem1} and, if $n$ is even, also the even valence 
property of lemma \ref{lem2}. A Gau\ss{} diagram of $c=2^n$ arrows,
with $G_n$ as intersection graph, would yield a value of $v_3$,
asymptotically equivalent modulo constants to $c\cdot (\log_2 c)^2$.
\end{exam}

Of course, a simple argument shows, that any graph containing
already $G_3$ as subgraph (i.~a., $G_4,G_6,\dots$) cannot be the
intersection graph of a Gau\ss{} diagram, but evidently we must invest more into the structure
of (intersection graphs of) Gau\ss{} diagrams. Unfortunately, the
further conditions will not be that simple and bringing them
into the game will make proofs (even more) tedious.

But, in any case, note, that the odd crossing number twist
knots ($!3_1,!5_2,!7_2,!9_2,\dots$) show, that we cannot prove more
than quadratical growth of $v_3$ in $c$. This is possibly, however,
indeed the worst case.

\begin{conj}
The positive twist knot diagrams minimize $v_3$ over all
connected irreducible positive diagrams of odd crossing number.
\end{conj}

\section{\label{un}Unknotting numbers and an extension
of Bennequin's inequality}

In this section, we introduce the machinery of Bennequin type
inequalities which will subsequently needed to prove further
properties of positive knots.

Parallelly we will also consider the following question,
which naturally arises in the study of knots (and links) via braids.

\begin{ques}\label{qu1}
Many classical properties of knots are defined by the
existence of diagrams with such properties. In how far
do these properties carry over, if we restrict ourselves to
closed braid diagrams?
\end{ques}

We will discuss this question in \S\ref{S.3} for positivity and,
applying our new criteria, give examples that the answer is in general
negative. On the other hand, here we will observe
for unknotting number the answer to be positive.


To start with, recall the result
of Vogel \cite{Vogel2} that each diagram is transformable into a 
closed braid diagram by crossing-augmenting Reidemeister
II moves on pairs of reversely oriented strands belonging to
\em{distinct} Seifert circles (henceforth called Vogel moves) only.

As observed together with T.~Fiedler, this result has 2 interesting
independent consequences. 
The first one is a ``singular'' Alexander theorem

\begin{theo}
Each $m$-singular knot is the closure of an  $m$-singular braid. 
\end{theo}

\proof Apply the Vogel algorithm to the $m$-singular diagram,
which clearly does not affect the singularities. \qed

This was, however, also known previously, see e.~g., \cite{Birman2}.

The other consequence is related to question \ref{qu1}.

\begin{theo}
Each knot realizes its unknotting number in a diagram as a 
closed braid. 

\end{theo}

\proof Take a diagram $D$ of $K$ realizing its unknotting number
and apply the Vogel algorithm obtaining a diagram $D'$.
As the crossing changes in $D$ commute with the Vogel moves,
the same crossing changes unknot $K$ in $D'$. \qed

I.~e., the answer of question \ref{qu1} for unknotting number
is yes! 

Combining Vogel's result with the Bennequin inequality \cite{Bennequin},
we immediately obtain

\begin{theo}\label{th4.3}
In each diagram $D$ of a knot $K$, 
\begin{eqn}\label{BVineq}
|w(D)|+1\,\le\, n(D)+2g(K)\,,
\end{eqn}
where $w(D)$ and $n(D)$ are the writhe and
Seifert circle number of $D$.
\end{theo}

\proof Bennequin proved the theorem for braid diagrams. From
this it follows for all diagrams by the Vogel algorithm, as
a Vogel move does not change neither the writhe nor the
number of Seifert circles. \qed

This fact for the unknot (which is also a special case of a result of
Morton \cite{Morton2}, who proved it for all achiral knots)
proves (in an independent way than theorem \ref{th1}) the following 

\begin{corr}
There is no non-trivial positive irreducible diagram
of the unknot.
\end{corr}

I.~e., in positive diagrams the unknot behaves as in
alternating ones.

\proof For such a diagram $D$, $n(D)=c(D)+1$, where $c(D)$ the
number of crossings of $D$.
Therefore, each smoothing of a
crossing in $D$ augments the number of components. Hence
no pair of crossings in the Gau\ss{} diagram can be linked, and so
all chords are isolated and all crossings are reducible. \qed

Another more general corollary is originally due to Cromwell
\cite{Cromwell}:

\begin{corr}\label{corr4.2}(Cromwell)
The Seifert algorithm applied to positive diagrams gives a
minimal surface.
\end{corr}

\proof This follows from theorem \ref{th4.3} together with
the formula for the genus of the Seifert algorithm surface
associated to $D$, which is $(c(D)-n(D)+1)/2$, as $|w(D)|=c(D)$ for $D$ positive.
\qed

Coming finally back to question \ref{qu1}, we see that we
have discussed the most interesting cases. For the crossing
number $10_8$ as an example as well. For braid
index the question does not make much sense, neither it does
for Seifert genus. Certainly the Seifert algorithm assigns
a surface to each diagram. However, Morton \cite{Morton2}
proved that there really exist knots, where in \em{no}
diagram the Seifert algorithm gives a minimal Seifert surface!
Posing question \reference{qu1} on minimality just for canonical
Seifert surfaces, that is, Seifert surfaces obtained by
the Seifert algorithm, the answer is again negative.
The knot $7_4$ has a positive diagram, and hence a
canonical Seifert surface of (minimal) genus $1$, whereas
by \cite{BoiWeb} the genus of a canonical Seifert surface in
any of its braid diagrams is minorated by its unknotting number $2$,
calculated by Lickorish \cite{Lickorish} and Kanenobu-Murakami
\cite{KanMur}.

The only interesting case to discuss is

\begin{ques}
Does each knot realize its bridge number in a closed braid diagram?
\end{ques}

Another question coming out of Vogel's result is

\begin{ques}
Does each knot realize its unknotting number in a diagram  as closed
braid of minimal strand number?
\end{ques}

%
%

Bennequin conjectured \eqref{BVineq} also to hold if we replace
genus by unknotting number (this is sometimes called the
Bennequin unknotting conjecture). This was recently proved by Kronheimer
and Mrowka \cite{KroMro} (see remarks below) and independently
announced by Menasco \cite{Menasco2}, but, to the best of my
knowledge, without a published proof. As before, Vogel's algorithm
extends this inequality.

\begin{theo}\label{th4.3a}
In each diagram $D$ of a knot $K$ it holds 
$|w(D)|+1\,\le\, n(D)+2u(K)\,.$
\qed\end{theo}

As we observed, \eqref{BVineq} is sharp for positive knots and
so we obtain

\begin{corr}\label{corr2}
For any positive knot $K$ it holds $u(K)\ge g(K)$. \qed
\end{corr}

This, combined with the inequality of Boileau-Weber-Rudolph
\cite{BoiWeb,Rudolph} leads to

\begin{corr}\label{corr3}
For any braid positive knot $K$ it holds $u(K)=g(K)$. \qed
\end{corr}

This was conjectured
by Milnor \cite{Milnor} for algebraic knots, neighborhoods of
singularities of complex algebraic curves, which are known to be
braid positive. 
%
%
Boileau and Weber \cite{BoiWeb} led it back to the
conjecture that the ribbon genus of an algebraic knot is equal
to its genus (see \S4 of \cite{Fiedler}), which was in turn known
by work of Rudolph \cite[p. 30 bottom]{Rudolph} to follow from the
Thom conjecture, recently proved by Kronheimer and Mrowka
\cite{KroMro}.

As pointed out by Thomas Fiedler, more generally, corollary
\reference{corr2} also follows
from Rudolph's recent result \cite{Rudolph3} that positive links
are (strongly) quasipositive, as he proved \cite{Rudolph2} that
a quasipositive knot bounds a complex algebraic curve in the 4-ball.
The genus of such a curve is equal to the lower bound for $g$ in
Bennequin's inequality and so not higher than $g$ itself.
Hence if a knot, which bounds a complex algebraic curve
is positive, then the the genus of the knot is equal to the
4-ball genus of the
complex algebraic curve that it bounds, which by \cite{KroMro}
was proved to realize the slice genus of the knot, and this is
as well-known always not greater than its unknotting number. 

Using corollary \reference{corr3} we can determine the
unknotting number of some knots.

\begin{exam}\label{_ex1}
The knots $10_{139}$ and $!10_{152}$ are braid positive, which is
evident from their diagrams in \cite{Rolfsen}. Their Alexander
polynomials tell us that they both have genus 4, hence their
unknotting number is also 4.
\end{exam}

Thus, we recover the result of Kawamura \cite{Kawamura}. However, 
corollary \reference{corr2} brings us a step further.

\begin{exam}\label{_ex2}
The knots $10_{154}$ and $10_{162}$ (the Perko duplication of $!10_
{161}$) are positive and have genus 3. Hence their unknotting number
is at least 3. Therefore, it is equal to 3, as 3 crossing changes
suffice to unknot both knots in their Rolfsen diagrams (find them!).
To determine the unknotting numbers in these examples is not possible
with the Bennequin unknotting conjecture for itself. Although both
knots satisfy \eqref{fib}, a property of braid positive knots we will
recall in \S\ref{S.3}, they are both not braid positive. As their genus
is 3, a positive $n$-braid realizing them would have $n+5$ crossings.
For $n<5$ this contradicts their crossing number, and for $n\ge 5$
such a braid would be reducible (getting us back to the case $n<5$).
\end{exam}

\begin{exam}\label{_ex3}
Another example is $!10_{145}$. $!10_{145}$ cannot be dealt with
directly by corollary \reference{corr2}, as it is not positive
(see \cite{Cromwell}). But it can be dealt with by observing that
is differs by one crossing change from $!10_{161}$, or by the original
Bennequin unknotting inequality: $!10_{145}$ is a (closed)
11 crossing 4-braid with writhe 7 \cite[appendix]{Jones2}.
Thus this knot  has unknotting number at least $2$.
On the other hand, 2 crossing changes suffice to unknot it as evident
from its Rolfsen diagram \cite[appendix]{Rolfsen}.
\end{exam}

For all 5 knots in examples \reference{_ex1}, \ref{_ex2} and \ref{_ex3}
the inequality $|\sg(K)/2|\le u(K)$ is not sharp, hence
the signature cannot be used to find out the unknotting number.
Therefore, this also disproves a conjecture of Milnor
(see \cite{Bennequin}), that $|\sg(K)/2|=u(K)$ for braid positive knots.
 
It is, however, striking that all 5 knots are non-alternating.
The reason for this is that if the positive diagram of $K$ is also
alternating, then indeed $|\sg(K)/2|=g(K)$, and hence (modulo
question \reference{qu_pa}) corollary \reference{corr2} (and even
the stronger corollary 1 of \cite{Rudolph3}) does not give anything
more for the unknotting number than the signature. One way to see
this is to use the principle of Murasugi and Traczyk (see
\cite{Traczyk} and \cite[p.~437]{Kauffman2}) to compute the signature
in alternating diagrams using the checkerboard shading and to observe
that if the alternating diagram is simultaneously positive, then
the white regions correspond precisely to the Seifert circles.

T.\ Kawamura informed me that some of the examples \reference{_ex1} and
\reference{_ex2} have been obtained independently by T.\ Tanaka
\cite{Tanaka}, who also found the unknotting number of $!10_{145}$,
inspiring me to give an independent argument in example
\reference{_ex3}. Very recently, A'Campo informed me that all
these examples have also been obtained independently by him
in \cite{ACampo}.

A further related, and meanwhile very appealing, conjecture was made in
\cite{MurPrz}. Using \cite[corollary 4.4]{restr}, it can be restated as
follows.

\begin{conj}
For any positive fibered knot $K$, $u(K)=g(K)$.
\end{conj}

We have seen this to be true for braid positive knots and also
for the two other positive fibered knots in Rolfsen's tables~--
$10_{154}$ and $!10_{161}$. The fact that a counterexample must have
$u>g$ makes the conjecture hard to disprove. Since for showing
$u>g$ any 4-genus estimate is useless, the only still handy
way would be to use Wendt's inequality \cite{Wendt} $t(K)\le u(K)$,
where $t(K)$ counts the torsion coefficients of the $\bZ$-homology
of the double branched cover of $S^3$ along $K$. We know from
the Seifert matrix that $t(K)\le 2g(K)$, and for example some 
(generalized) pretzel knots show, that this inequality is sometimes
sharp. Thus knots
with $t(K)>g(K)$ exist. However, they are very special and indeed
there was no positive fibered prime knot of $\le 16$ crossings
with $t(K)>g(K)$ (even $t(K)=g(K)$, where the methods of \cite{unkn1,%
Traczyk2} may have a chance to work, if all torsion coefficients are
divisible by $3$ or by $5$, was satisfied only by the trefoil).

\section{Further properties of the Fiedler Gau\ss{} sum invariant}

Here are two properties of $v_3$ which we will conclude with.

\begin{theo}
If $K$ is a positive knot, then $v_3(K)\ge 4g(K)$.
\end{theo}

\proof Take a positive diagram of $K$. As both the genus
of the canonical Seifert surface (which we observed in \S\ref{un}
is minimal for positive diagrams) and $v_3$ are additive
under connected sum of diagrams, assume that the diagram is
non-composite. Furthermore assume w.l.o.g., that the diagram cannot
be reduced by a Reidemeister I move after eventually
previously performing a sequence of Reidemeister III moves,
so it is in particular bireduced (else reduce the
diagram this way, noting that by the above remark this procedure
does not change the genus of the canonical Seifert surface).

So we can assume, we have a non-composite bireduced positive
diagram of $c$ crossings and $n$ Seifert circles. If $n=1$ the
diagram is an unknot diagram and the result is evident.
If $n=2$ the diagram is of a $(2,m)$-torus knot $K_{2,m}$,
$m$ odd and the result follows from a direct calculation of $v_3$
on $K_{2,m}$ (noting that $g(K_{2,m})=\frac{|m|-1}{2}$).
So now assume $n\ge 3$. Then the genus of the Seifert surface is
\[
g(K)=\frac{c-n+1}{2}\,\le\,\frac{c-2}{2}
\]
Therefore $4g(K)\le 2c-4$. But on the other hand by theorem \ref{th1},
$v_3(K)\ge 2c-4$. \qed

\begin{theo}
Let $D$ be a positive reduced 
diagram and $D'$ be obtained from $D$
by change of some non-empty set of crossings. Then
$v_3(D')<v_3(D)$.
\end{theo}

This fact may not be too surprising, as $v_3$ in general
increases with the number of positive crossings. However, it
is not obvious in view of the fact, that $v_3$ sometimes decreases
when a negative crossing is switched to a positive one.

\proof To compare $v_3(D')$ and $ v_3(D)$, we need to figure
out how the configurations in \eqref{v3} change by switching
the positive crossings in $D$.

The configurations of the first two terms in \eqref{v3} remain
in $D'$ (as orientation of the arrows does not matter) but possibly
change their weight. In any case the weight of such a
configuration in $D'$ is not higher than (the old  weight)
$1$ and so the contribution of these two
configurations to the value of $v_3$ decreases from $D$ to $D'$.

Something more interesting happens with the third term. 
A configuration in $D$ may or may not survive in $D'$.
But even if it does, its weight in $D'$ is not more than one.
However, a new configuration of positive weight can be created
in $D'$. It happens if it has exactly two negative arrows
and they are linked. This we will call an interesting configuration.
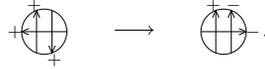
\[
\conf1142653+++\quad\lra\quad\conf1416253--+\,.
\]

To deal with the interesting configurations, we will find other ones
whose negative contributions equilibrate these of interesting
configurations. First note, that any interesting configuration
has a canonical pair of a negative arrow $p$ and a half-arc $c$
assigned:
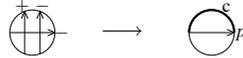
\[
\conf1416253--+\quad\lra\quad
\GD{\labar{180}{0}{p}
\piclinewidth{10}\piccirclearc{0 0}{1}{0 180}
\picputtext{1.3 60 polar}{c}
}
\]
Now consider any such canonical pair in the GD together with
all arrows starting outside and ending on the half-arc $c$.
Assume there are $l$ negative and $m$ positive such arrows.

Now beside the interesting configurations there are several ones
which bring a decrease of $v$ (see figure \ref{fig1};
$p$ is always the arrow from left to right and $c$ the upper half-arc).


\let\@nl\\
\begin{figure}[htb]
\begin{tabular}{*4c}
configuration in $D'$ & times appearing & 
difference of contributions & times counted \\
& &  to $v_3$ from $D$ to $D'$ & \\
 \conf1416253--+ \& \conf1416253-+- & $\le l \cdot m$ & $+1$ & 1\\
\conf1416253--- \& \conf0415263--- & together $\mybin{l}{2}$
 & $-2$ & 1, 2\\
\conf1416253-++ \& \conf0415263-++ & together $\mybin{m}{2}$
 & $-2$ & 1\\
\conf33142-- & $l$ & $-2$ & 2\\
\conf33142-+ &  $m$ & $-1$ & 1
\end{tabular}
\caption{\label{fig1}}
\end{figure}

To compute the total contribution of all these configurations
to the change of value of $v_3$ for one specific canonical pair,
we have to multiply their number with the difference of contributions,
dividing by the number of counting them with respect to different
 canonical pairs. The resulting contribution for the configurations 
in figure \ref{fig1} is for a given canonical pair
\[
\le l\cdot m-\mybin{l}{2}-2\mybin{m}{2}-l-m\,=\,
-\frac{(m-l)^2}{2}-\frac{l}{2}-\frac{m^2}{2}\,,
\]
which is negative for $l,m\ge 0$ unless $l=m=0$. But $p$
in a canonical pair with $l=m=0$ is reducible. This shows the
theorem. \qed

A classical result on the Jones polynomial \cite{Kauffman3,%
Murasugi,Thistle} states that a non-composite alternating
and a non-alternating diagram of the same crossing number
never belong to the same knot. This is no longer true,
if we replace `alternating' by `positive', as we will observe in
\S\ref{qu}. However, it is true if instead of non-compositeness
we demand the diagrams to have the same plane curve.

\begin{corr}
The Jones polynomial always distinguishes a reduced positive
and a non-positive diagram with the same plane curve. \qed
\end{corr}

\section{\label{secCas}The Casson invariant on positive knots}

Here we shall say a word on the degree-2-Vassiliev
invariant $v_2$, sometimes attributed to Casson because of its relation
to the 3-manifold invariant discovered by him, see \cite{Akbulut}.
This invariant is the coefficient of $z^2$ in the Conway polynomial
$\nabla(z)$, or alternatively $\Delta''(1)/2$, where
$\Delta$ is the Alexander polynomial \cite{Alexander}.
Using the Polyak-Viro formula for it,
we obtain a similar result for positive knots as for $v_3$.

\begin{theo}\label{theo1}
In a positive reduced $c$ crossing diagram, $v_2\ge c/4$.
\end{theo}

This bound is sharp, as again the four crossing diagram of the
(positive) trefoil shows. However, under assumption of
bireducedness more seems possible.

\begin{conj}\label{cj3}
In a positive bireduced diagram $D$, $v_2\ge lk(D)/4$, where $lk(D)$
is the number of linked pairs in $D$.
\end{conj}

As a consequence, if this conjecture is true, for example
in a positive bireduced $c$ crossing diagram, $v_2\ge c/3$.
The reason why this bound is suggestive lies in the method we introduce
to prove theorem \reference{theo1} and it will be motivated later.

The proof of theorem \ref{theo1} uses the Polyak-Viro formula for $v_2$
\begin{eqn}\label{eq2}
v_2\,=\,\frac 12\left(\GD{\arrow{225}{45}
\arrow{315}{135}
\picfillgraycol{0}\picfilledcircle{1 90 polar}{0.09}{}
} + \GD{\arrow{45}{225}
\arrow{135}{315}
\picfillgraycol{0}\picfilledcircle{1 90 polar}{0.09}{}
} \right)\,,
\end{eqn}
obtained by symmetrization from the formula \cite[(3)]{VirPol}.

A similar but somewhat more complicated formula for
$v_2$ was found by Fiedler \cite{Fiedler,Fiedler4}, who uses it to show
that for a braid positive knot $K$ it holds $v_2(K)\ge g(K)$. This implies
theorem \ref{theo1} for braid positive knots because of the inequality
\begin{eqn}\label{eq1}
g(K)\ge c/4
\end{eqn}
in a reduced braid positive $c$ crossing diagram $D$ of a braid positive
knot $K$, which is a consequence of both the Bennequin inequality
\cite[theorem 3, p.~101]{Bennequin} and Cromwell's work on homogeneous
knots \cite[corollaries 5.1 and 5.4]{Cromwell}, see \cite{genera}
or remark \reference{r3.5}.

Note, that \eqref{eq1} is not true in general for positive knots,
but Fiedler's inequality extends to this case.

\begin{theo}\label{th4}
For positive knots it holds $v_2(K)\ge g(K)$ and $v_2(K)\ge u(K)$.
\end{theo}

\begin{rem}
By corollary \reference{corr2}, the first inequality in theorem 
\reference{th4} follows from the second one. But the argument of our
proof shows it without involving the slice version of Bennequin's
inequality. Therefore, we felt it deserves in independent exposition.
\end{rem}
Note, that this excludes a large class of positive (see \cite{Cromwell})
polynomials as Conway polynomials of positive knots,
$1+z^2+z^4$ is a simple one (belonging, inter alia, to the knot
$6_3$).

In the sequel, we will need the following fact, which we invite the
reader to prove.

\begin{exer}\label{ex0}
Show that in the Gau\ss{} diagram of a positive knot diagram any arrow is
distinguished in exactly half of the pairs in which it is linked.
\end{exer}

\proof[of theorem \ref{theo1}] Call a linked pair in a based Gau\ss{}
diagram \em{admissible}, if it is of one of the two kinds appearing
in \eqref{eq2}.

Fix a reduced positive diagram $D$. We apply now a series of 
transformations to $D$ we call \em{loop moves}, ending at the
trivial diagram. What is crucial for our argument, is that
\eqref{eq2} is independent of the choice of a base point. That means,
as long as we can assure that the Gau\ss{} diagram is realizable,
that is, corresponds to a knot diagram, we can place for the next
loop move on the diagram the base point to some other favorable place.

Now we describe how to perform a loop move. Take a crossing $p$ in $D$
whose smoothing produces a (diagram of a link with a) component
$K$ with no self-crossings. In the Gau\ss{} diagram this means,
that the arrow of $p$ does not have non-linked arrows on both its
sides. In $D$, $p$ looks like this
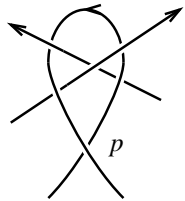
\[
\diag{1cm}{3}{3.5}{
  \pictranslate{-1 0}{
    \picPSgraphics{0 setlinejoin}
    \picline{1.5 1.5}{2.25 2}
    \picline{3.5 1.8}{2.2 2.45}
    \picmulticurve{0.12 1 -1.0 0}{2 0.5}{2.5 1}{3 2}{3 2.3}
    \picmulticurve{0.12 1 -1.0 0}{3 0.5}{2.5 1}{2 2}{2 2.3}
    \picellipsevecarc{2.5 2.3}{0.5 0.7}{0 90}
    \picellipsearc{2.5 2.3}{0.5 0.7}{90 180}
    \picmultivecline{0.12 1 -1.0 0}{2.2 2.45}{1.5 2.8}
    \picmultivecline{0.12 1 -1.0 0}{2.25 2}{3.75 3}
    \picputtext[l]{2.8 1.2}{$p$}
  }
}\,.
\]
Now switch appropriately at most (but, in fact, exactly) half of the
crossings on the loop and remove them from the diagram:
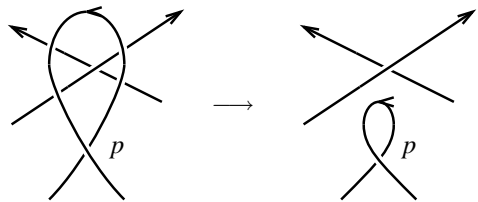
\[
\diag{1cm}{3}{3.5}{
  \pictranslate{-1 0}{
    \picPSgraphics{0 setlinejoin}
    \picline{1.5 1.5}{2.25 2}
    \picline{3.5 1.8}{2.2 2.45}
    \picmulticurve{0.12 1 -1.0 0}{2 0.5}{2.5 1}{3 2}{3 2.3}
    \picmulticurve{0.12 1 -1.0 0}{3 0.5}{2.5 1}{2 2}{2 2.3}
    \picmultivecline{0.12 1 -1.0 0}{2.2 2.45}{1.5 2.8}
    \picmultivecline{0.12 1 -1.0 0}{2.25 2}{3.75 3}
    \picmultiellipsevecarc{0.12 1 -1.0 0}{2.5 2.3}{0.5 0.7}{0 90}
    \picmultiellipsearc{0.12 1 -1.0 0}{2.5 2.3}{0.5 0.7}{90 180}
    \picputtext[l]{2.8 1.2}{$p$}
  }
}\,\lra\,
\diag{1cm}{3}{3.5}{
  \pictranslate{-1 0}{
    \picPSgraphics{0 setlinejoin}
    \picline{1.5 1.5}{2.25 2}
    \picline{3.5 1.8}{2.2 2.45}
    \picmulticurve{0.12 1 -1.0 0}{2 0.5}{2.5 1}{2.7 1.2}{2.7 1.5}
    \picmulticurve{0.12 1 -1.0 0}{3 0.5}{2.5 1}{2.3 1.2}{2.3 1.5}
    \picellipsevecarc{2.5 1.5}{0.2 0.3}{0 90}
    \picellipsearc{2.5 1.5}{0.2 0.3}{90 180}
    \picmultivecline{0.12 1 -1.0 0}{2.2 2.45}{1.5 2.8}
    \picmultivecline{0.12 1 -1.0 0}{2.25 2}{3.75 3}
    \picputtext[l]{2.8 1.2}{$p$}
  }
}\,.
\]
In the Gau\ss{} diagram this corresponds to removing the $k$ arrows
linked with $p$. Note, that by even valence $2|k$ and by reducedness
$k>0$. Assume in the resulting diagram $D'$ there are $c$ reducible
crossings, $p$ including. Any of them must have been linked in $D$
with at least 2 crossings linked with $p$ (because of even valence).

Now, we place the basepoint in the Gau\ss{} diagram as follows:
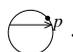
\[
\GD{\labar{160}{20}{p}
\picfillgraycol{0}\picfilledcircle{1 40 polar}{0.09}{}
}\,.
\]
Henceforth, such a picture means that there is no other end of an arrow
between the basepoint and the arrow end to which it is depicted to be
close (here the over-crossing of $p$).

Remove the arrows linked with $p$ in $D$
by exercise \ref{ex0} removes $k/2$ admissible linked pairs
with $p$ and for any other reducible crossing $p'$ in $D'$,
at least one admissible linked pair ($p'$ must have been linked
in $D$ with some even non-zero number of to be removed arrows and
by exercise \reference{ex0} exactly half of them gives with
it an admissible linked pair). Then the
procedure of building $D'$ out of $D$ and reducing $D'$ reduces the
value of $v_2$ at least by
\[
\frac 12\left(\frac k2+c-1\right)\,,
\]
an hence by integrality of $v_2$ at least by
\[
v_2(D)-v_2(D')\,\ge\,\br{\frac k4 + \frac c2}\,,
\]
whereas it reduces the crossing number of the (reduced)
diagram by $k+c$. The ratio
\[
\frac {k+c}{\br{\frac k4 + \frac c2}}
\]
for $2|k$, $k,c>0$ is at most $4$, unless $k=4$ and $c=1$ (in
which case it is $5$). We would like to show that in this case
$v_2$ reduces at least by 
\begin{eqn}\label{k2c}
\frac{1}{2}\left(\frac k2+c\right)\,,
\end{eqn}
and hence by integrality at least by $\ds \br{\frac k4 + \frac {c+1}2}$.

To do so, now consider some $p''\ne p$, which is not linked with $p$
and does not become reducible in $D'$. The loop move reduced the
number of arrows linked with $p''$ by some even number $2l$, possibly
$0$, such that half of this number (that is, $l$) of arrows point
in either direction w.r.t. $p''$. Then for each such $p''$ the loop move
reduces $v_2$ additionally by $l/2$. What we need is that at least for one
crossing $p''$ in $D$ we have $l>0$. This occurs, unless $p$ belongs
to a connected component of $D$, in which all $p''\ne p$ are linked
with $p$ or linked with all $p'$ linked with $p$. The connected
component would have $c+k$ crossings and would be resolved by the loop
move (and the elimination of reducible crossings following it).

But in our case $c+k=5$ and on the two positive diagrams of 5 crossings
$v_2$ is $2$ (for $5_2$) resp.\ $3$ (for $5_1$). Therefore,
\eqref{k2c} follows.

Resolving this case, we have always ensured $v_2(D)-v_2(D')\ge (k+c)/4$,
and so the theorem follows inductively over $c(D)$, as it is true
for $c(D)=0$ and any positive diagram can be trivialized by a sequence
of the above transformations. \qed

\begin{corr}
For any of the polynomials of Alexander/Conway, Jones, HOMFLY,
the Brandt-Lickorish-Millett-Ho polynomial $Q$ \cite{BLM,Ho}
and Kauffman \cite{Kauffman2} only finitely many positive knots
have the same polynomial and there is no positive knot with
unit polynomial.
\end{corr}

\proof Use the inequality for $v_2$ and the relations
\[
-6v_2\,:=-3\Delta''(1)=V''(1)=Q'(-2)
\]
and the well-known specializations for the HOMFLY and Kauffman
polynomial. The equality between the Jones and Alexander polynomial is
probably due already to Jones \cite[\S 12]{Jones2}. The relation
between the Jones and Brandt-Lickorish-Millett-Ho polynomial is
proved by Kanenobu in \cite{Kanenobu3}. \qed

\begin{rem}
As a consequence of the result of \cite{Kauffman3,Murasugi,Thistle}
on the span of the Jones polynomial, only finitely many alternating
knots have the same Jones polynomial. On the other hand, indeed
collections of such knots (sharing even the same HOMFLY and Kauffman
polynomials) of any finite size exist \cite{Kanenobu}. It would be
interesting whether similar constructions to these of Kanenobu are
also possible
in the positive case for both the Jones and Conway/Alexander polynomial
and also to give an infinite series of \em{alternating} knots having the
same Conway/Alexander polynomial, similar to the one (of non-alternating
knots) in \cite{Kanenobu2}. Note, that knots of such a series (except
finitely many) can neither be skein equivalent nor (by \cite{Cromwell}) 
fibered. 
\end{rem}

In any case, Kanenobu's examples of \cite{Kanenobu} show that the
lower bound for the crossing number coming from the span of the
Jones polynomial can be arbitrarily bad. Theorem \ref{theo1} gives
us a new tool for positive knots.

\begin{corr}
Let $K$ be a knot with a positive reduced diagram of $c$ crossings.
Then $c(K)\ge \sqrt{2c}$.
\end{corr}

\proof Use \cite[theorem 2.2.E]{VirPol2}. \qed

Although we will sharpen it, we already remark the inequality
$v_2(K)\ge g(K)/2$ we obtain for the genus of a positive knot $K$ from
the inequality $g\le c/2$.


\begin{exer}\label{ex1.5}
Prove that if $D$ is a positive reduced diagram and $D'$ is obtained
from $D$ by change of some but not all of its crossings, then
$v_2(D')<v_2(D)$.
\end{exer}

\proof[of theorem \ref{th4}]
We use again the inductive step in the proof of theorem
\ref{theo1}. Fix a loop in a positive diagram $D$ bounded by a crossing
$a$. Assume the loop has $2c$ crossings on it. Then, switching at most
(but, in fact, exactly) $c$ crossings on the loop, it can be pulled
above or below all the strands intersecting it.

Now recall the inequality of Bennequin-Vogel \eqref{BVineq} of
\S\reference{un}.
The inequality is sharp for $D$ positive, as shows the (therefore
minimal) surface coming from the Seifert algorithm. This shows, that
the switching of the $c$ crossings in $D$ reduces the absolute writhe
maximally by $2c$, and so (as it does not affect $n(D)$) $g$ at
\em{most} by $c$. On the other hand, as we will just
observe, it reduces $v_2$ at \em{least} by $c$. The following
Reidemeister moves do not change $v_2$ or $g$ and
then the same inductive argument as in the proof of theorem \ref{theo1} 
applies to show the first inequality asserted in the theorem.
For the second one note, that the procedure describes an unknotting
of $K$ (and hence the number of crossing changes is at least its
unknotting number).

To see that removing the arrow of $a$ and all its linked arrows
in the Gau\ss{} diagram to $D$ reduces $v_2$ at least by $c$,
put the basepoint in the Gau\ss{} diagram as follows:
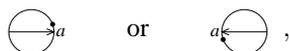
\[
\GD{\labar{180}{0}{a}
\picfillgraycol{0}\picfilledcircle{1 20 polar}{0.09}{}
}\qquad\text{ or }\qquad
\GD{\labar{0}{180}{a}
\picfillgraycol{0}\picfilledcircle{1 200 polar}{0.09}{}
}\,,
\]
and use the Polyak-Viro formula
\begin{eqn}\label{PV1}
v_2\,=\,\GD{\arrow{225}{45}
\arrow{315}{135}
\picfillgraycol{0}\picfilledcircle{1 90 polar}{0.09}{}
}
\end{eqn}
together with exercise \ref{ex0}. \qed

\begin{exer}\label{ex2}
Modify the proof of theorem \ref{th4}
to show that in a positive diagram $D$, $lk(D)\ge 3g$, and deduce from this the
inequality for any arbitrary diagram.
\end{exer}

\hint Use that beside the linked pairs with $a$
any arrow linked with $a$ must be linked with
another arrow in the Gau\ss{} diagram to $D$.

We finish the discussion of $v_2$ in its own right by an inequality
involving both the crossing and unknotting number of a positive diagram.

\begin{theo}\label{theocu}
Let $D$ be a reduced positive diagram of crossing number $c(D)$
and unknotting number $u(D)$. Then $\ds v_2(D)\ge \frac{c(D)+u(D)}{5}$.
\end{theo}

\begin{rem}
Replacing the `5' in the denominator by `4', theorem \reference{theocu}
would imply theorem \reference{theo1}, and replacing the `5' by `6',
it would follow from it using $u(D)\le c(D)/2$. Thus `5' is in a sense
indeed the interesting denominator. On the other hand, for braid
positive knots theorem \reference{theo1} indeed follows from
theorem \reference{theocu} because of theorem \reference{thm4}, a
property of braid positive knots we are going to prove later.
\end{rem}

\proof We split the proof into two steps recorded as several lemmas.
Our strategy will be as follows.
{\def\labelenumi{\theenumi.}\mbox{}\\[-18pt]

\begin{enumerate}
\item\label{step1} Apply loop moves to $D$, that do not unknot any
connected component of $D$, until you obtain a diagram $D''$ with the
property that any of its connected components gets unknotted by any
loop move on it. Show the inequality of theorem \reference{theocu}
for $D''$.
\item\label{step2} Show that if $D'$ arises from $D$ by a loop move
of step \reference{step1}, then 
\[
5\bigl(\,v_2(D)-v_2(D')\,\bigr)\,\ge\,c(D)-c(D')+\myfrac{k}{2}\,,
\]
where $\myfrac{k}{2}$ is the number of crossings switched by the loop
move (so $k$ is the number of crossings on the loop). The totality
of all such crossings switches over all moves of step \reference{step1}
together with any unknotting sequence for $D''$ forms an
unknotting operation of $D$ (because the removal of a loop after a loop
move commutes with all subsequent crossing changes), and the length
of this unknotting sequence is $\ge u(D)$.
\end{enumerate}
}

All diagrams we consider in the sequel will be assumed positive.
For the first step we need some preparation.

\begin{lemma}\label{l1}
Let $D$ be a connected diagram on which \em{any} loop move unknots.
We call $D$ loop-minimal. Then and only then $D$ (more exactly its
Gau\ss{} diagram) has no subdiagram of the kind $(5,1)$:
$ \GD{\chrd{-60}{60}
\chrd{-130}{90}
\chrd{-90}{130}} $.
\end{lemma}

\proof Call an arrow corresponding to a crossing on which a loop move
can be applied extreme. If $D$ has a configuration of the kind $(5,1)$
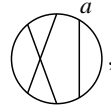
\[
\CD{\labch{-60}{60}{a}
\chrd{-130}{90}
\chrd{-90}{130}}\,,
\]
then we can w.l.o.g. find $a$ to be extreme and applying a loop move
to $a$ we get a diagram with a linked pair, which is hence knotted.
This contradiction shows the direction `$\So$'. To show `$\Longleftarrow
$', observe that the miss of a configuration of the kind $(5,1)$ in $D$
means that if for two arrows $p$ and $p''$ in $D$, $p\ncap p''$, then
it holds $\fa p':\,p'\cap p''\So p'\cap p$.

But then all $p''$ remain reducible after a loop move on $p$.
This shows the other direction, letting $p$ vary over all arrows of $D$.
\qed

\begin{lemma}\label{l2}
Let $D$ be prime and loop-minimal. Then there exists a (necessarily
disjoint) decomposition $\{\mbox{ arrows of $D$ }\}=K\cup L$, such that
$K\ne\varnothing\ne L$ and 
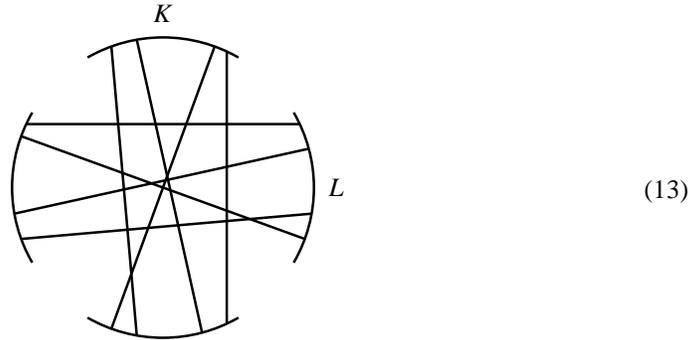
\begin{eqn}\label{eqnKL}
\fa p\in K,\,p'\in L: p\cap p'\,.
\end{eqn}
\begin{eqn}\label{picKL}
\diag{1cm}{5}{5}{
\pictranslate{2.56 2.56}{
  \picmultigraphics[rt]{4}{90}{
    \piccirclearc{0 0}{2}{-30 30}
  }
  \picputtext[d]{0 2.2}{$K$}
  \picputtext[l]{2.2 0}{$L$}
  \picmultigraphics[rt]{2}{90}{
    \picline{-65 2 x polar}{65 2 x polar}
    \picline{-110 2 x polar}{70 2 x polar}
    \picline{-75 2 x polar}{100 2 x polar}
    \picline{110 2 x polar}{-100 2 x polar}
  }
}
}
\end{eqn}
Moreover, unless $D$ is the 3 crossing trefoil diagram, $K$ and $L$
can be chosen to have $\ge 2$ elements each.
\end{lemma}

\proof Distinguish two cases.
\begin{caselist}

\case There is no triple of arrows of type $(3,3)$. Then by
straightforward arguments $D$ has the form
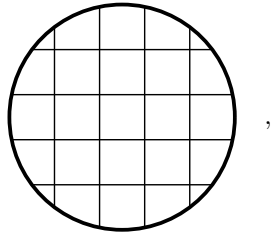
\[
\diag{1.5cm}{2}{2}{
  \pictranslate{1 1}{
    \piccircle{0 0}{1}{}
    \piclinewidth{80}
    \picclip{\piccircle{0 0}{1}{}}
       {\picmultigraphics[rt]{2}{90}{
          \picmultigraphics{4}{0.4 0}{\picline{-0.6 -1}{-0.6 1}}
        }
       }
  }
}\quad,
\]
that has the desired property. By even valence $|K|$ and $|L|$ are both
even, so $\ge 2$.

\case There is a triple $(a,b,c)$ of arrows of type $(3,3)$.
Number the segments of the base line the 6 ends of the arrows
separate by $1$ to $6$:\\[2mm]
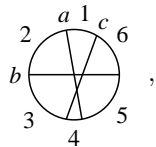
\[
\CD{
 \labch{-100}{60}{c}
 \labch{0}{180}{b}
 \labch{-80}{100}{a}
 \picputtext{1.4 40 polar}{6}
 \picputtext{1.4 80 polar}{1}
 \picputtext{1.4 140 polar}{2}
 \picputtext{1.4 225 polar}{3}
 \picputtext{1.4 270 polar}{4}
 \picputtext{1.4 320 polar}{5}
}\quad,
\]\\[2mm]
and let for $i,j\in\{1,\dots,6\}$
\[
[i,j]\,:=\,\{\mbox{ chords with endpoints on segments $i$ and $j$ }\}\,.
\]
By lemma \reference{l1}, $[i,j]=\varnothing$ if $|i-j|\in\{0,1,5\}$,
hence there are 9 possibilities left for $i,j$. Then we can write
down some $K$ and $L$ ad hoc. Set
\[
K\,:=\,\{\,b,\,c,\,[1,3],\,[3,5],\,[3,6],\,[4,6],\,[2,6]\,\}
\quad\mbox{and}\quad
L\,:=\,\{\,a,\,[1,4],\,[1,5],\,[2,4],\,[2,5]\,\}\,,
\]
and verify the desired property case by case. For example, any
element $p\in [1,4]$ must intersect any element $p'$ in $[1,5]$,
for if not, a loop move in $p'$ would preserve the linked pair $(a,p)$.

We may have now that $|L|=1$. But whenever we have a decomposition with
\eqref{eqnKL} and $K\ne\varnothing\ne L$, we can build a new one by
taking some $p\in L$ and setting $K':=\{\,p'\,:\,p'\cap p\,\}$ and $L':=
\{\,p'\,:p'\ncap p\,\land \,\fa p''\,:\,p''\cap p\iff p''\cap p'\,\}$.
Then $K'$ and $L'$ still satisfy \eqref{eqnKL} and by even valence
$|K'|\ge 2$. If now $|L'|=1$, that is, $L'=\{p\}$ for any choice of $p$,
then any two arrows are linked in $D$, and it is a $(2,2n+1)$ torus
knot diagram. But for such a diagram $K'$ and $L'$ with $|K'|,|L'|\ge 2$
are immediately found, unless $n=1$, which is the 3 crossing trefoil
diagram. 
\end{caselist}
\qed

\begin{exer}\label{exf}
Show that in fact a Gau\ss{} diagram with out a configuration of type
$(5,1)$ are either rational knot diagrams of the form $C(p,q)$ with
$p$ and $q$ even integers, or (generalized) pretzel diagrams
$P(a_1,\dots,a_n)$, $n$, $a_i$ odd.
\end{exer}

\hint Consider a chord $p_1$ with the maximal valence 
(number of linked chords) and collect $p_1$ and all its non-linked
chords into a collection.  Then consider from the rest of the
chords again one chord $p_2$ with maximal valence and so on.
You obtain a decomposition of the chords into collections, such that
two chords intersect if and only if they belong to distinct
collections. Use even valence to show that either there is an odd
number of collections each one of odd size, which is the pretzel
diagram case, or an even number of collections each one of even size,
and in this case deduce that the are no more that two collections
using the non-realizability of 
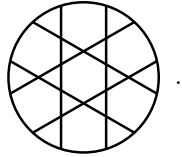
\[
\diag{1cm}{2}{2}{
  \pictranslate{1 1}{
    \piccircle{0 0}{1}{}
    \picclip{\piccircle{0 0}{1}{}}
       {\picmultigraphics[rt]{3}{60}{
	  \picmultigraphics{2}{0.6 0}{\picline{-0.3 -1}{-0.3 1}}
        }
       }
  }
} \enspace.
\]

\begin{lemma}\label{l3}
If $D$ is connected, reduced and loop-minimal, then $v_2(D)\ge \ds
\frac{c(D)-2}{2}$\,.
\end{lemma}

\proof Assume $c(D)\ge 4$, as the trefoil diagram is easily checked.
Then by lemma \reference{l2} we have the decomposition of the arrows
into $K$ and $L$. Then in the picture \eqref{picKL} we can w.l.o.g.,
modulo rotating the diagram by $90^\circ$ (swopping $K$ and $L$)
and mirroring, assume that $\ge \myfrac{k}{2}$ arrows in $K$
point upward, and $\ge \myfrac{l}{2}$ arrows in $L$ point from left
to right. Then placing the basepoint above the arrows in $L$
and to the right of the arrows in $K$ and using the formula
for $v_2$ in \eqref{PV1}, we see
\[
v_2\,\ge\,\frac{kl}{4}\,\ge\,\frac{k(c-k)}{4}\,\ge\,
\frac{2(c-2)}{4}\,=\,\frac{c-2}{2}\,,
\]
as $k\ge 2$, $c-k\ge 2$ and $c\ge 4$. \qed

\begin{lemma}\label{l4}
(first step) If $D$ is a loop-minimal connected prime or composite
diagram, then $v_2(D)\,\ge\,\ds\frac{c(D)+u(D)}{5}$\,.
\end{lemma}

\proof Assume first $D$ is prime. As $u(D)\le c(D)/2$ it suffices to 
show $v_2\ge\,\frac{3c(D)}{10}$. But lemma \reference{l3} gives
$v_2\ge\myfrac{c(D)}{2}-1$, which is better, unless $c(D)<5$,
which is directly checked.

If $D$ is composite, use that $c(D)$, $u(D)$ and $v_2(D)$ are additive
under connected sum of \em{diagrams} (bewaring that the question for 
the crossing and unknotting number for \em{knots} is a 100 year
old conjecture, that no one knows how to prove except in special
cases!).  \qed

\begin{lemma}\label{l5}
(second step) If a loop move $D\to D'$ does not unknot a connected
component of $D$, then
\[
5\bigl(\,v_2(D)-v_2(D')\,\bigr)\,\ge\,c(D)-c(D')+\myfrac{k}{2}\,,
\]
where $\myfrac{k}{2}$ is the number of crossings switched by the loop
move (so $k$ is the number of crossings on the loop). 
\end{lemma}

\proof By the proof of theorem \reference{th4} we have
\[
v_2(D)-v_2(D')\,\ge\,\frac k2\,,
\]
and by the proof of theorem \reference{theo1} we have
\[
v_2(D)-v_2(D')\,\ge\,\frac k4+\frac c2\,,
\]
again denoting by $c$ the number of reducible crossings after the move.
Hence by arithmetic mean
\begin{eqnarray*}
v_2(D)-v_2(D') & \ge & \frac 12\,\left(\frac k4+\frac c2+\frac k2\right)\\
	       & = & \frac {3k}{8} + \frac c4 \\
	       & = & \left(\frac {c+k}{4}\right) + \frac {k/2}{4}\,,
\end{eqnarray*}
and as $c+k=c(D)-c(D')$, already
\[
4\bigl(\,v_2(D)-v_2(D')\,\bigr)\,\ge\,c(D)-c(D')+\myfrac{k}{2}\,,
\]
that certainly remains true when replacing the factor `4' by `5'. \qed

Using the strategy outlined in the beginning, lemmas \reference{l4}
and \reference{l5} prove theorem \reference{theocu}. \qed

\begin{rem}
It is striking that the whole proof goes through with denominator
`4' instead of `5', except at one point: the case $c(D)=4$ in
lemma \reference{l4} (the positive 4 crossing trefoil diagram). This
is, however, for our argument fatal, because we would need to control
how many such factors occur in the diagram $D''$ after step \reference{%
step1}. One hope to get out of the dilemma would be to find
loop moves, such that connectedness is always preserved, but 
one can find examples, where this is not possible. 

Moreover, along similar lines one shows that $3v_2(D)$
decreases not slower than $c(D)$ under the moves of step
\reference{step2}. So the motivation for conjecture \reference{cj3}
is again the problem how to handle step \reference{step1}.
Similarly to connectedness, it is difficult to make the loop
move behave well w.r.t. bireducedness.
\end{rem}

\begin{corr}
If $K$ is positive, then $5v_2(K)\ge\max\deg V(K)$.
\end{corr}

\proof Use that by \cite{Kauffman3,Murasugi,Thistle}
$\sp V(K)\le c(D)$ on a positive (or any other) diagram $D$ of $K$, and
that $u(D)\ge u(K)\ge g(K)=\min\deg V(K)$ by corollary 
\reference{corr2} and \cite[theorem 4.1]{restr}. \qed

It is interesting to remark that this is a entirely
combinatorial statement that heavily relies on this deep topological
fact~-- the truth of the (local) Thom conjecture. It would
be nice to know whether it cannot be derived also completely
combinatorially.

\begin{rem}
The work done in this paragraph in an easy manner also recovers
for positive knots the mentioned result of Cochran--Gompf and Traczyk
on the positivity of the signature. For this it suffices to remark
that a loop move, consisting of switching positive crossings
to negative, never reduces the signature, and that it is
positive on the knots of exercise \ref{exf} by direct calculation.
\end{rem}

\section{Relations between $v_2$, $v_3$ and the HOMFLY polynomial%
\label{rp}}

The Polyak-Viro-Fiedler formulas also allow to relate both
the degree-2 and degree-3 Vassiliev invariants to each other
in positive diagrams, giving a lower bound for their crossing number.

In the following we give inequalities resulting from 
such combined applications of the various formulas, which,
while not terribly sharp, hardly seem provable using other arguments.

Let
\[
l_i\,:=\,\#\,\{\text{ crossings linked with crossing $i$ }\}
\]
in some fixed positive diagram $D$ of $c$ crossings.

\begin{lemma}
In a positive diagram $D$ of $c$ crossings,
\begin{eqn}\label{eq3}
v_3\,\ge\,\sum_{i=1}^c\,\mybinn{l_i/2}{2}+\frac{l_i}{2}\,=\,
\sum_{i=1}^c\,\frac{(l_i+1)^2-1}{8}\,.
\end{eqn}
\end{lemma}

\proof We have in a positive diagram 
\[
v_3(D)\,\ge\,
\GD{\arrow{225}{45}
\arrow{315}{135}
}\,+\,
\conf1{-4}15362{i}{}{}
\,+\,
\conf0415263{i}{}{}\,.
\]

The first term on the right gives the second summand
in \eqref{eq3} (note, that a linked pair is counted twice for both
arrows in it). Numbering the horizontal chord in terms
2 and 3 by $i$, we see that the sum of terms
2 and 3 is the count of pairs of equally oriented arrows with respect
to arrow $i$. Now by exercise \ref{ex0} for each $i$ and for each
orientation there are two collections of $l_i/2$ equally
oriented arrows with respect to arrow $i$, giving $2\mybinn{l_i/2}{2}$
possible choices of pairs of equally oriented arrows. As the cases
where the equally oriented
arrows are linked are counted twice, we factor out the `2' and obtain
the formula \eqref{eq3}. \qed

\begin{lemma}\label{lemma2}
In a positive diagram $D$,
\begin{eqn}\label{eq4}
4 v_2\,\le\,\sum_{i=1}^c\,{l_i}\,.
\end{eqn}
\end{lemma}

\proof This is obviously a consequence of \eqref{eq2}. \qed

\begin{theo}\label{theo2}
In a positive reduced $c$ crossing diagram, $c>0$, then $v_3>v_2$ and
\[
c\,\ge\,\frac{2v_2^2}{v_3-v_2}\,.
\]
\end{theo}

\proof In view of \eqref{eq4}, the right hand side of \eqref{eq3}
is minimized by $l_i:=4v_2/c$, in which case it becomes
\[
\frac{c\ds\left(\frac{4v_2}{c}+1\right)^2}{8}-\frac{c}{8}\,.
\]
So
\begin{eqn}\label{eq5}
v_3\,\ge\,\frac{2v_2^2}{c}+{v_2}\,,
\end{eqn}
from which the assertion follows, as by
\eqref{eq2} and \eqref{v3} always $v_3>v_2$ (even $v_3\ge 2v_2$). \qed

\begin{rem}
The Polyak-Viro formula for $v_3/4$
\[
\frac{v_3}{4}\,=\,\frac 12 \conf1415326{}{}{} + \conf0145236{}{}{}
\]
proves similarly as \eqref{eq3} the inequality
\[
\frac{v_3}{4}\,\le\,\frac 12 
\sum_{i=1}^c\,\left(\frac{l_i}{2}\right)^2\,.
\]
Adding $\ds{\frac{v_3}{4}}-\sum\frac{l_i^2}{8}$ to
the right hand side of \eqref{eq3}, we obtain
$\ds\frac 34 v_3\,\ge\,\sum_{i=1}^c\,\frac{l_i}{4}$, so
by lemma \ref{lemma2}, $\ds v_3\,\ge\,\frac{4v_2}{3}$. This
inequality is weaker than \eqref{eq5}, if $v_2\,\ge\,\frac{c}{6}$,
which we showed always holds in reduced positive diagrams.
\end{rem}

\begin{theo}\label{theo3}
In a positive diagram of $c$ crossings,
\[
\frac 34\,v_3\,\le\,v_2\,\cdot\,c\,.
\]
\end{theo}

\proof As before, combining the Fiedler and Polyak-Viro formulas,
we have
\[
\frac 34\,v_3\,=\,
\GD{\arrow{225}{45}
\arrow{315}{135}
}\,+\,
\conf1{-4}15362{}{}{}
\,+\,
\GD{\arrow{225}{45}
\arrow{315}{135}
\arrow{270}{90}}
\,-\,\frac 12\conf1415326{}{}{}\,\le\,
\GD{\arrow{225}{45}
\arrow{315}{135}
}\,+\,
\conf1415362{}{}{}\,+\,
\conf1145362{}{}{}\,+\,
\GD{\arrow{225}{45}
\arrow{315}{135}
\arrow{270}{90}}\,.
\]
But for positive diagrams the four terms on the right are
equal to the first four terms in the numerator of the $v_2$ formula
\cite[p.~5 bottom]{VirPol}. \qed

All previous calculations suggest that in the point of view
of Gau\ss{} sums, the following invariant plays some key role.

\begin{defi}
Define the linked pair number $lk(K)$ of a knot $K$ as the minimal
number of linked pairs in all its diagrams.
\end{defi}

What can we say about $lk(K)$? As a consequence of \eqref{eq2}
or \eqref{eq4}, $lk(K)\ge 2v_2(K)$ for any knot $K$. However,
$v_2(K)$ may be sometimes negative. A non-negative lower bound
is $3g(K)$, following from exercise \reference{ex2}. In fact, exercise
\reference{ex2} shows $lk(K)\ge 3\tl g(K)$, where $\tl g(K)$
is the weak Seifert genus of $K$, that is, the minimal
genus of a surface, obtained by applying the Seifert algorithm
to any diagram of $K$. As we noted, sometimes $\tl g(K)>g(K)$.
Morton showed \cite{Morton2}, that
\[
\tl g(K)\ge \max\deg_{m}\,P(K)/2\,,
\]
so
\[
lk(K)\ge \frac 32\,\max\deg_{m}\,P(K)\,.
\]
On the other hand, for $K$ positive we proved
$v_3(K)\ge \frac 43 lk(K)$, so we obtain a self-contained inequality

\begin{prop}
For a positive knot $K$ we have
\begin{myeqn}{\qed}
v_3(K)\ge 2\max\deg_{m}\,P(K)\,.
\end{myeqn}
\end{prop}

This condition is also violated by our previous example $12_{2038}$.
We also obtain

\begin{prop} $\ds v_3(K)\ge \frac 83v_2(K)$ for $K$ positive. \qed
\end{prop}

As simple examples show, except for the low crossing number cases
and connected sums thereof these inequalities are far from being sharp,
so significant improvement seems possible. The problem with pushing
further our inductive arguments in \S\ref{secCas} is that
it appears hard to control how often these
low crossing number cases occur as connected components in intermediate
steps of trivializing a positive diagram with our move.

A final nice relation between $v_2$ and $v_3$ is unrelated to Gau\ss{}
sums and bases on an idea of Lin. Let $w_{\pm}$ denote the untwisted
double operation of knots with positive (resp. negative) clusp.

\begin{prop}
$v_3(w_{\pm}(K))=\pm 8 v_2(K)$.
\end{prop}

\proof The dualization $w^*_\pm$ of $w_\pm$ is a nilpotent
endomorphism of $\cV^n$, the space of Vassiliev invariants of
degree at most $n$. But $\cV^3/\cV^2$ and $\cV^2/\cV^1$ are
one-dimensional and hence are killed by $w^*_\pm$. Therefore,
$w^*_\pm$ maps $v_2$ to a constant and checking it on the unknot
we find that it is zero (this also follows from $\Delta=1$ for
an untwisted Whitehead double of any knot). $v_3$ is taken to
something in degree at most 2, so $v_3(w_\pm(K))=c_1^\pm\,v_2(K)+
c_0^\pm$. $c_0^\pm=0$ follows from taking the unknot and to see
$c_1^\pm=\pm 8$ check that $v_3$ is 8 an the positive clusp
untwisted Whitehead double of one of the trefoils. \qed

Combining this with our Gau\ss{} sum inequalities we immediately
obtain

\begin{corr}
An untwisted Whitehead double of a positive knot has non-self-conjugate
Jones polynomial. In particular, the knot is chiral and
has non-trivial Jones polynomial. Moreover, there are only finitely many
positive knots, whose untwisted Whitehead doubles (or similarly,
twisted Whitehead doubles with any fixed framing) have the same
Jones polynomial. \qed
\end{corr}

\section{Braid positive knots\label{S.3}}

The following section deals with the more specific subclass of
positive knots, namely those with positive braid representations.
First, as a digression from the Gau\ss{} sum approach,
we improve some inequalities of Fiedler \cite{Fiedler} on the
degree of the Jones polynomial of such knots, and latter we write down
certain inequalities for the Casson invariant of knots with
positive braid representations, giving some applications.

{\bf Notation.} For a braid $\beta$ denote by $\hat\beta$ 
its closure, by $n(\beta)$ ist strand number and by $[\beta]$
its homology class (or exponent sum), i.~e. its image under the
homomorphism $[\,.\,]\,: B_n\to H_1(B_n)=B_n/\langle [B_n,B_n] \rangle
\simeq \bZ$, given by $[\sg_i]=1$, where $\sg_i$ are the Artin
generators.

\begin{defi}\label{defipa}
A knot is called braid positive, if it has a positive diagram as
a closed braid.
\end{defi}

\note The term ``braid positive'' is self-invented and provided to
give a naturally seeming name for such knots and links, distinguishing
them from the ones we call `positive'. However, braid positive
knots are called sometimes ``positive knots'' elsewhere in
the literature, so beware of confusion!


First we will recall and sharpen an obstruction of Fiedler
\cite{Fiedler} to braid positivity.

\begin{@lemma}[\cite{Fiedler}]\rm
For any braid positive $k$ component link $L$ without trivial split
components, we have $\min\deg V(L)>0$ and $\min\cf V(L)=(-1)^{k-1}$.
\end{@lemma}

Here is our improved version of Fiedler's result.

\begin{theo}\label{thm4}
If $L$ is a non-split $k$ component link, $L=\hat\beta$, with $\beta$ a
positive reduced braid of $c$ crossings, then 
\begin{eqn}\label{**}
\min\deg V(L)\ge c/4-\frac{k-1}{2}\,\ge\, c(L)/4-\frac{k-1}{2}
\end{eqn}
and $\min\cf V(L)=(-1)^{k-1}$.
\end{theo}

To prove the theorem, let's start with the

\begin{figure}[htb]
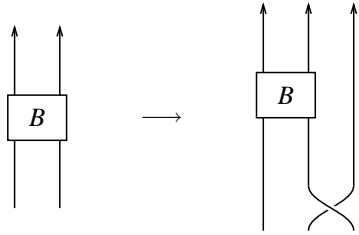

\[
\diag{6mm}{3}{4}{
  \picvecline{1 0}{1 4}
  \picvecline{2 0}{2 4}
  \picfilledbox{1.5 2}{1.3 1}{$B$}
}\quad\lra\quad
\diag{6mm}{4}{5}{
  \rbraid{2.5 0.5}{1 1}
  \picvecline{1 0}{1 5}
  \picvecline{2 1}{2 5}
  \picvecline{3 1}{3 5}
  \picfilledbox{1.5 3}{1.3 1}{$B$}
}
\]
\caption{\label{MarII}The Markov II move.}
\end{figure}

\begin{lemma}\label{lem3}
If a positive braid diagram of a prime knot is reducible,
then it admits a reducing Markov II \cite{Birman} move, see figure
\reference{MarII}.
So, if a prime knot has a positive (closed) braid diagram,
it also has a reduced one.
\end{lemma}

\proof Take a reducible crossing in the closed braid diagram
and smooth it out. As the knot is
prime, assume w.l.o.g. that the right one of the two resulting closed
braid diagrams belongs to the unknot. If we know, that each
positive braid diagram of the unknot is either
trivial or reducible, repeat this procedure, ending up with a
trivial (braid) diagram of the unknot on the right. Then the
last smoothed crossing is one corresponding to a reducing
Markov II move.

For positive braids it follows from work of Birman and
Menasco \cite{BirMen} and also from the Bennequin inequality
\cite[theorem 3, p.~101]{Bennequin}, that if $\hat\beta$ is the unknot,
then $|[\beta]|<n(\beta)$. Therefore, if $\beta$ is positive,
it must contain each generator exactly once, so all its
crossings are reducible. \qed

\begin{rem}\label{rm1}
Note, that our capability to control so well positive braid diagrams
of the unknot by these (deeper) results, is rather surprising, as
in general there exist \em{extremely ugly} braid diagrams of the unknot
\cite{Morton,Fiedler2}.
\end{rem}

\begin{rem}
A similar statement is also true for alternating diagrams.
To see the fact, that each alternating braid diagram of the unknot is
either trivial or reducible, recall the result
of Kauffman \cite{Kauffman3}, Murasugi \cite{Murasugi} and
Thistlethwaite \cite{Thistle}, that all alternating
diagrams of the unknot are either trivial or reducible.
\end{rem}

The assertion in lemma \reference{lem3} in the positive
case is also true for composite knots and links.

\begin{lemma}\label{lem4}
Any braid positive link has a reduced braid positive diagram.
\end{lemma}

\proof
In the braid positive diagram use the iteration of the procedure
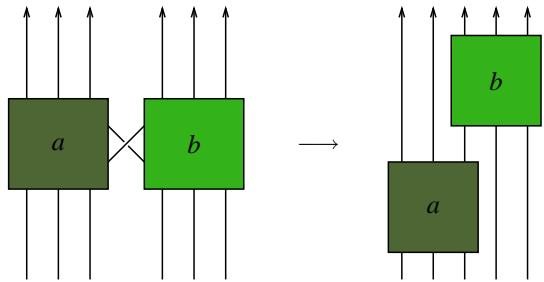
\[
\diag{6mm}{6}{6}{
  \picmultigraphics{2}{3 0}{
    \picmultigraphics{3}{0.7 0}{
      \picvecline{0.8 0}{0.8 6}
  }}
  \picline{3.4 2.6}{2.6 3.4}
  \picmultiline{0.12 1 -1 0}{2.6 2.6}{3.4 3.4}
  {\picfillgraycol{0.3 0.4 0.2}
   \picfilledbox{1.5 3}{2.2 2}{$a$}
  }
  {\picfillgraycol{0.2 0.7 0.1}
   \picfilledbox{4.5 3}{2.2 2}{$b$}
  }
} \quad\lra\quad
\diag{6mm}{4}{6}{
  \picmultigraphics{5}{0.7 0}{
    \picvecline{0.6 0}{0.6 6}
  }
  {\picfillgraycol{0.3 0.4 0.2}
   \picfilledbox{1.3 1.6}{2 2}{$a$}
  }
  \pictranslate{2.7 4.4}{
    \picscale{1 1}{
      \picfillgraycol{0.2 0.7 0.1}
      \picfilledbox{0 0}{2 2}{$b$}
  }}
}
\]
which gives a reduced diagram. \qed
%

\begin{rem}
Note, that, however, for alternating diagrams the above procedure
does not work. The granny knot $!3_1\#!3_1$ has a
reducible alternating diagram as
closed 4-braid, but no alternating diagram as closed 3-braid.
\end{rem}

\proof[of theorem \ref{thm4}]
Take equation (10) of \cite{Fiedler} for positive $\beta$.
\begin{eqn}\label{(7)}
\min\deg V(L)=\frac{1}{2}([\beta]+1-n(\beta))
\end{eqn}
As $\beta$ is w.l.o.g. by lemma \reference{lem4}
reduced, and generators appearing only once in $\beta$
correspond to reducible crossings in the closed braid diagram, we have 
\begin{eqn}\label{***}
[\beta]\ge 2(n(\beta)-k),
\end{eqn}
so
\begin{eqn}\label{*}
\min\deg V(L)\ge \frac{n(\beta)-k}{2}\,.
\end{eqn}

On the other hand, as $\beta$ positive, $[\beta]=c$, so
\[
\min\deg V(L)\ge \frac{c+1-n(\beta)}{2}\,.
\]
Therefore
\[
\min\deg V(L)\ge \min_n \max\left( \frac{n(\beta)-1}{2}\,,\,
\frac{c+1-n(\beta)}{2} \right)-\frac{k-1}{2}\,=\,\frac{c}{4}-\frac{k-1}{2}\,.
\]
The second assertion follows directly from \cite[theorem 2]{Fiedler}. \qed

\begin{rem}
Applying $n\ge b(\hat\beta)$ in \eqref{*}, or taking 
the inequality $c(L)\ge 2(b(L)-k)$ of Ohyama \cite{Ohyama}
in \eqref{**}, we also obtain the weaker
\[
\min\deg V(\hat\beta)\ge \frac{b(\hat\beta)-k}{2}\,.
\]
\end{rem}

\begin{rem} Considering $L=K$ to be a knot,
the first inequality in \eqref{**} is evidently sharp, as a braid with
each generator appearing twice shows. Concerning the
second inequality and demanding the braid to be irreducible (i.~e. not
conjugate to a braid with an isolated generator),
the inequality \eqref{***} can be further improved a little
by observing, that a positive braid with exactly $2(n(\beta)-1)$
crossings is still transformable modulo
Yang-Baxter relation (that is, a transformation of the kind
$\sg_{i-1}\sg_{i}\sg_{i-1}=\sg_{i}\sg_{i-1}\sg_{i}$)
into one with isolated generators.
So we can add a certain constant on the r.h.s. of \eqref{***},
and to our bound, maybe excluding some low crossing cases
($!3_1$ and $!5_1$ show that $\BR{c(K)/4}$ at least is sharp.)

However, at $c(K)/4+2$ there will be really something to do, as
for $[\beta]=2n(\beta)+6$ there is a series of examples of braids
$\{\beta_n\,|\,n \text{ odd}\,\}$ with 
\[
\beta_n=\bigl((\sg_1\sg_3\dots\sg_{n-4}\sg_{n-2}^3)\,
	(\sg_2^3\sg_4\dots\sg_{n-3}\sg_{n-1})\bigr)^2
\]
or schematically
\[
\def\bh#1{\hbox to 1ex{\hss$#1$\hss}}
\beta_n\,=\,\begin{array}{*9c}
  & 3 &   & 1 & \cdots &   & 1 &   & 1 \\
1 &   & 1 &   & \cdots & 1 &   & 3 &   \\
  & 3 &   & 1 & \cdots &   & 1 &   & 1 \\
1 &   & 1 &   & \cdots & 1 &   & 3 &   \\[2mm]
\bh{\sg_1} & \multicolumn{7}{c}{\leaders\hbox to 1ex{\hss.\hss}\hfill} &
\bh{\sg_{n-1}}
\end{array}
\]
which do not admit a Yang-Baxter relation modulo cyclic permutation
and close to a knot. Of course, this is far away from saying
that $\beta_n$ are irreducible or that even $\hat\beta_n$ is a minimal
diagram (which would mean, that the second bound is also sharp)
but I don't know how to decide this.
\end{rem}

\begin{rem}\label{r3.5}
The expression appearing on the r.h.s. of \eqref{(7)}
is equal to
\[
g(L)+\frac{1-n}{2}
\]
where $n$ is the number of components of $L$. This follows
from the (classical) formula for the genus of the canonical
Seifert surface, together with the fact that this
Seifert surface is (of) minimal (genus)
in positive diagrams, see corollary \reference{corr4.2}.
(This observation is treated in more detail and generalized in
\cite{restr}). Therefore, for a braid positive knot $K$ we have
\begin{eqn}\label{fib}
\max\deg\Delta(K)\,=\,g(K)\,=\,\min\deg V(K)\,\ge\, c(K)/4\,,
\end{eqn}
where $\Delta$ is the Alexander polynomial and the first
equality comes from the fiberedness of the knot. The condition
\eqref{fib} is not sufficient, though. We have seen this in
example \reference{_ex2}.
\end{rem}

As a braid positive knot $K$ by lemma \ref{lem3}
always has a reduced braid positive diagram, and
a reduced braid positive diagram by theorem \ref{thm4} does not more
than $4\min\deg V(K)$ crossings, we see that  braid positivity
can always be decided. This, of course, works with the results of
the previous section as well, but this bound is \em{considerably}
sharper.

Here we shall observe that braid positive
is really stronger than positive, so our definition
\ref{defipa} is justified.

\begin{exam}
The knot $!5_2$ is positive (see, e.~g., \cite{Kauffman2}).
We have (see, e.g., \cite[Appendix]{Adams} or \cite{polynomial})
$\min\deg V(!5_2)=1<\myfrac{5}{4}$. So, although positive,
$!5_2$ can never be represented as a closed positive braid.
The same is true for $!7_2$ and $7_4$. Note, that in all 3 examples the
conclusion of non-braid positivity would not have been possible
with Fiedler's weaker criterium.
\end{exam}

\begin{rem}
It is known, that closed positive braids are fibered, and that
fibered knots have monic Alexander polynomial \cite[p.~259]{Rolfsen}
(i.~e., with edge coefficients $\pm 1$), so the monicness of the
Alexander polynomial is also an obstruction to braid positivity,
and applies in the above 3 examples $!5_2,!7_2$ and $7_4$ as well.
Another way to deal with these cases is to use the observation, that
they all have genus 1 (which can be seen by applying the Seifert
algorithm to their alternating diagrams \cite{Gabai}), and the fact
(following from the Bennequin inequality
\cite[theorem 3, p.~101]{Bennequin}), that the only braid positive
genus 1 knot is the positive trefoil. A special way to exclude $7_4$
is to use that it has unknotting number $2$ \cite{Adams}, contradicting
the inequality $u(K)\le g(K)$ for braid positive knots $K$
due to Boileau and Weber \cite{BoiWeb} and Rudolph
\cite[prop. on p. 30]{Rudolph}, see also \cite{Bennequin}.
\end{rem}

\begin{exam}
The 10 crossing knot $!10_2$ is fibered and his minimal degree of
the Jones polynomial is positive, but it is 1, so $!10_2$ is not
a closed positive braid.
$!10_2$, however, can also be dealt with by the non-positivity
of its Conway polynomial \cite{Busk}.
\end{exam}

\begin{exam}
On the other hand, the knots $7_3$ and $!7_5$ are positive, but
their minimal Jones polynomial degree $2$ does not
tell us, that they are not braid positive. But they have non-monic
Alexander polynomial, and so they cannot even be fibered.
\end{exam}

The variety of existing obstructions to (braid) positivity makes it hard
to find a case, where our condition is universally better. Here is
a somewhat stronger example, coming out of some quest 
in Thistlethwaite's tables.

\begin{figure}[htb]
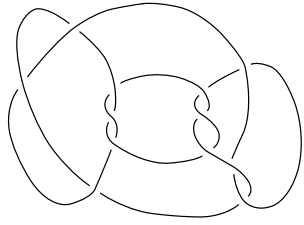

\[
\epsfs{4cm}{k-12-1930}
\]
\caption{\label{fig12-1930}The knot $!12_{1930}$.}
\end{figure}

\begin{exam}
The knot $!12_{1930}$ on figure \reference{fig12-1930} has the
HOMFLY polynomial
\[
(4l^8+2l^{10}-l^{12})+(-4l^4+2l^6-4l^8+l^{10})m^2+l^4m^4\,.
\]
It shows, that no one of the above mentioned (braid) positivity
obstructions of \cite{Rolfsen,Busk,%
Cromwell,MorCro,Fiedler} is violated, but ours is. However, although
monic, the Alexander polynomial can be indirectly used to
show non-braid positivity. How?
\end{exam}

We now give some improvements of the inequalities for positive knots
for closed positive braids of given strand number. It is obvious
that without this restriction not more than a linear lower bound for
$v_2$ and $v_3$ in $c$ can be expected, as shows the iterated connected
sum of trefoils (we will shortly construct more such examples).

\begin{theo}\label{thv2bp}
If $\bt$ is a positive braid of exponent sum (or crossing number)
$[\bt]$, and $n$ strands, closing to a knot, then
\begin{eqn}\label{poi}
v_2(\hat \bt)\,\ge\,\frac{[\bt]^2}{4n(n-1)}-\frac{(2n-3)(n-1)}{8}\,.
\end{eqn}
\end{theo}

\proof Consider $l_{ij}=lk(i,j)$ for $1\le i<j\le n$, the linking number
of strands $i$ and $j$ in $\bt$. Then $[\bt]=\sum_{i<j}l_{ij}$ and
each pair of strands $i$ and $j$ contributes to the Gau\ss{} diagram of
$\hat\bt$ a collection of $l_{ij}$ mutually linked arrows. If $l_{ij}$
is odd, the contribution of these arrows to the Gau\ss{} sum is
(independently of the choice of basepoint) the one of the $(2,l_{ij})$
torus knot, namely $(l_{ij}^2-1)/8$, while for $l_{ij}$ even, the
contribution is dependent of the choice of basepoint (changes by
$\mp 1$), but is in any case at least $(l_{ij}^2-4)/8$. The bound
on the right of \eqref{poi} is obtained by taking all $l_{ij}$ equal,
namely $[\bt]\big/{n\choose 2}$,
and using that at least $n-1$ of the $l_{ij}$ are odd, as any one-%
cycle permutation of $n$ elements has length at least $n-1$. \qed

\begin{corr}
If $\bt$ is a braid of $n$ strands, closing to a knot, 
$[\bt]_-$ the number of negative crossings in $\bt$, and
$[\bt]_0=[\bt]+2[\bt]_-$ the total number of crossings of $\bt$, then
\[
v_2(\hat \bt)\,\ge\,\frac{[\bt]_0^2}{4n(n-1)}-\frac{(2n-3)(n-1)}{8}-
\frac{[\bt]_-\bigl([\bt]_0-[\bt]_-\bigr)}{2}\,.
\]
\end{corr}

\proof Use the expression of $v_2$ in \eqref{eq2}, showing that
switching $[\bt]_-$ positive crossings in any diagram of $[\bt]_0$
crossings, decreases $v_2$ at most by the third term on the right. \qed

This means, that for $[\bt]_-$ sufficiently small, we have $v_2(\hat
\bt)>0$, implying as before that in particular $\hat\bt$ has
non-trivial $\Dl$, $V$ and $Q$ polynomial, and untwisted Whitehead
doubles with non-trivial $V$ polynomial.

In a similar way one proves

\begin{theo}\label{thv3bp}
If $\bt$ is a positive braid of $n$ strands, closing to a knot, then
\begin{eqn}
v_3(\hat \bt)\,\ge\,C_1\frac{[\bt]^3}{n^4}-C_2n^2\,,
\end{eqn}
for some (effectively computable and independent on $\bt$ and $n$)
constants $C_{1,2}>0$. \qed
\end{theo}

Theorems \ref{thv3bp} and \ref{thv2bp} imply a positive solution to
Willerton's problem 5 in \cite[\S 4]{Willerton} for positive
braids of given strand number.

\begin{corr}\label{cr32}
If $(\bt_i)$ are distinct positive braids of $n$ strands, then
$v_3(\hat\bt_i)\asymp v_2(\hat\bt_i)^{3/2}$, in particular,
$\lim\limits_{i\to\infty}\log_{v_2(\hat\bt_i)}v_3(\hat\bt_i)=\myfrac{3}
{2}$. \qed
\end{corr}

Such a property can be used to show that certain special positive braids
are not Markov equivalent to positive braids of given strand number,
where the calculation of the Homfly polynomial (and the
bound of the Morton-Williams-Franks inequality \cite{Morton2,WilFr})
can be tedious.

\begin{exam}
Let $\bt_i$ and $\bt_i'$ be positive braids of length $O(i^{1/4-\eps})$
in some $B_n$ with $\phi(\bt_i)=(1\,2)$, $\phi\,:\,B_n\to S_n$ being the
permutation homomorphism. Let $\{\,.\,\}_j$ be the shift map
$\sg_i\mapsto \sg_{i+j}$. Set 
\[
\bt_{[j]}=\prod\limits_{i=1}^j \{\bt_i\}_{i-1}\,\cdot\,\{\bt_j'\}_j\,
\in\,B_{n+j}\,,
\]
so that $\phi(\bt_{[j]})=(1\,2\,\dots\,n+j)^{-1}$. The all but
finitely many of the knots $\hat\bt_{[j]}$ have no positive
braid representations of some fixed strand number. To
see this, use that any crossing in $\hat\bt_{[j]}$ has
$O(j^{1/4-\eps})$ linked crossings, from which
the $v_3$ formula shows $v_3(\hat\bt_{[j]})=O(j^{3/2-2\eps})$, so
if $\log_{v_2(\hat\bt_{[j]})}v_3(\hat\bt_{[j]})\to x$, we must have
$x\le \myfrac{3}{2}-2\eps$, contradicting corollary \ref{cr32}.
\end{exam}

This example also shows that Willerton's problem cannot be solved
positively in general for braid positive knots. We can already take
the iterated connected sum of trefoils, but we also see how to
construct prime examples using \cite{Cromwell2}. For example,
take all $\bt_i=\sg_1^2\sg_2^2\sg_1\in B_3$ and $\bt_i'=\sg_1$,
and you get $v_2(\hat\bt_{[j]})=O(j)=v_3(\hat\bt_{[j]})$.

\section{\label{qu}Questions on positive knots}

After alternating knots have been well understood, it's interesting
to look for another class of knots. The positive knots provide
many interesting questions in analogy to alternating knots.

Here are some appealing questions thinking on alternating knots.

By \cite{Kauffman3,Murasugi,Thistle} any alternating reduced diagram
is minimal. We saw that, ignoring the second reduction move, this is not
true for positive knots. Is it true with the second reduction move
(and all its cablings)?
It seems, however, that things are not that easy with
positive knots (or the other way round~-- it makes them the more
challenging!).

\begin{exam}\label{ex8}
Consider the knot, which is the closed rational tangle with the
Conway notation $(-1,$ $-2,-1,-2,-5)$. Its diagram as
closed $(-1,-2,-1,-2,-5)$ tangle is
reduced and alternating, and hence minimal. The knot, however, 
has also a positive diagram as closed $(1,2,1,2,1,1,-1,-3)$ tangle,
which is bireduced, but non-minimal.
(This is one of a series of such examples I found by a
small computer program.)
\end{exam}



Conversely, for alternating prime knots, any minimal diagram
is alternating. As the example of the Perko pair \cite[fig.~10]
{Kauffman2} shows, this is not true for positive knots.
So we can ask:

\begin{ques}\label{qu5.1}
Does any positive knot have at least \em{one}
positive minimal diagram? If so, is there a set of
local moves reducing a positive diagram to a positive minimal diagram?
\end{ques}

%
%

\begin{rem}
I tried to find counterexamples to question \ref{qu5.1}
using the following (common) idea:
Consider the Conway notation $a=(a_1,\dots,a_n)$ of a (diagram of a)
rational tangle $A$, closing to a positive (diagram of some) knot
$K$. Then take some expression $c=(c_1,\dots,c_m)$ of its iterated
fraction
\[
a_n+\frac{1}{a_{n-1}+
	     \ds\frac{1}{\ds a_{n-2}+
		      \frac{1}{a_{n-3}+\dots
		              }
                     }
	    }
\,=\,
c_m+\frac{1}{c_{m-1}+
	    \ds \frac{1}{\ds c_{m-2}+
		      \frac{1}{c_{m-3}+\dots
		              }
                     }
            }
\]
with all $c_i$ of the same sign. The (diagram of the) tangle $C$
with Conway notation $c$ is equivalent to $A$ \cite{Adams}, closes to
an alternating diagram $\bar C$ of $K$. $K$ is also prime
(e.~g. by \cite{Menasco}, as $\bar C$ is non-composite and alternating).
Therefore any minimal diagram of this knot is alternating and
if $\bar C$ is not a positive diagram, by Thistlethwaite's invariance
of the writhe \cite{Kauffman2} it would follow that, as $K$
has one non-positive minimal diagram, no minimal diagram can be
positive (and also it didn't matter which $c$ you chose).
My computer program revealed, however, that there is no such
$a$ with $|a|\le 26$ (where $|a|:=\sum_{i=1}^n\,|a_i|$; note, that
by minimality of alternating diagrams always $|c|<|a|$).
Is there such an $a$ at all?
\end{rem}

\begin{ques}
Is (something like) the Tait flyping conjecture \cite{MenThis}
true for positive knots, i.~e. are minimal positive diagrams
transformable by flypes?
\end{ques}

\begin{ques}\label{qu5.3}
Menasco \cite{Menasco,Adams}~/~Aumann \cite[p.~150]{Adams} proved that
composite/split alternating links appear composite/split in \em{any}
alternating diagram. Using the linking number, it's easy to see that
for split links latter is also true in the positive case. But what is
with composite knots?
\end{ques}

In view of corollary \ref{corr4.2}, this is
a special case of a conjecture of Cromwell \cite[conjecture 1.6]
{Cromwell2}. Note, that affirming questions \ref{qu5.1} and \ref{qu5.3}
we would prove the additivity of the crossing number for positive
knots under connected sum.

\begin{ques}
If question \ref{qu5.3} has a negative answer, is still the weaker
statement true that positive composite knots have (only or at least one)
positive prime factor(s)?
\end{ques}

\begin{ques}\label{qu_pa}
Is it possible to classify alternating positive knots? 
Does an alternating positive knot always have 
a (simultaneously) alternating (and) positive diagram?
(Note, that this question for prime knots and question \ref{qu5.1} for
prime alternating knots are the same.)
\end{ques}

A question on unknotting numbers is

\begin{ques}
Does any positive knot realize its unknotting number
in a positive diagram?
\end{ques}

If the answer were yes, by arguments  analogous to those in the proof
of theorem \reference{th4}, the inequality of Bennequin-Vogel
\eqref{BVineq} would show that $u(K)\ge g(K)$ independently
from Menasco's result, so it is consistent with it.
%
%
%

A final question is suggested by the comparison between the
growth rates of $v_2$ and $v_3$ on positive knots.

\begin{ques}
What can be said about the sets
\[
S:=\{\,\log_{v_2(K)}v_3(K)\,:\,K\ne !3_1\,\mbox{positive}\,\}
\]
and
\[
SB:=\{\,\log_{v_2(K)}v_3(K)\,:\,K\ne !3_1\,\mbox{braid positive}\,\}\,?
\]
We have shown that $1\in\tl S=\bar S\sm\disc S\subset[1,3]$, and
similarly for $\widetilde{SB}$, where $\bar S$ denotes closure and
$\disc S$ the discrete subset of points of $S$. Is $\tl S\subset [1,2]$
or even $\tl S=[1,2]$? Is $\widetilde{SB}$ equal to or at
least contained in $[1,\myfrac{3}{2}]$?
\end{ques}

\noindent{\bf Acknowledgement.} I would wish to thank to W.\ B.\ R.\ 
Lickorish and J.~Birman for their kindness and help and especially to
T.~Fiedler for his remarks and pointing out an extension of theorem
\ref{th4.3}. I am grateful to T.\ Kawamura for informing me about
the results of T.\ Tanaka and to T.\ Kanenobu for informing me about
\cite{Kanenobu3}.

{\small

}

\end{document}